\documentclass[11pt,fleqn,twoside]{article}
\usepackage{amsfonts,amssymb,amscd,latexsym,slava-99}

\newtheorem{theorem}{Theorem}[section]
\newtheorem{remark}{Remark}[section]
\newtheorem{proposition}{Proposition}[section]
\newtheorem{conjecture}{Conjecture}[section]
\newtheorem{example}{Example}[section]
\newtheorem{definition}{Definition}[section]
\newcommand{\bea}{\begin{eqnarray}}
\newcommand{\eea}{\end{eqnarray}}
\newcommand{\beano}{\begin{eqnarray*}}
\newcommand{\eeano}{\end{eqnarray*}}
\newcommand{\bdm}{\begin{displaymath}}
\newcommand{\edm}{\end{displaymath}}
\newcommand{\pf}{{\it Proof:\ }}
\newcommand{\epf}{\hfill \fbox{}}

\begin{document}

\thispagestyle{empty}

\renewcommand{\ehkol}{O.S. Stoyanov}
\renewcommand{\ohkol}{Poisson Diffeomorphism Groups}

%\newtheorem{theorem}{Theorem}
%\newtheorem{proposition}{Proposition}
%\newtheorem{lemma}{Lemma}

%\begin{flushleft}
%\footnotesize \sf
%Journal of Nonlinear Mathematical Physics \qquad 1999, V.6, N~3,
%\pageref{kupershmidt_5-fp}--\pageref{kupershmidt_5-lp}.
%\hfill {\sc Article}
%\end{flushleft}

\vspace{-5mm}

\renewcommand{\footnoterule}{}
{\renewcommand{\thefootnote}{}
 \footnote{\prava{O.S. Stoyanov}}}

\name{Poisson Diffeomorphism Groups}\label{kupershmidt_5-fp}

\Author{Ognyan S. STOYANOV}

\Adress{~~P.O.Box 2833,  Champaign, IL 61825, USA\\
~~E-mail: stoyanov@biobig.rutgers.edu}

\Date{~~April 5, 2000}

\begin{abstract}
\noindent
We construct explicitly a class of coboundary Poisson-Lie structures on the group of formal diffeomorphisms of ${\Bbb R}^n$.
Equivalently, these give rise to a class of coboundary triangular Lie bialgebra structures on the Lie algebra $W_n$ of formal vector
fields on ${\Bbb R}^n$. We conjecture that this class accounts for all such coboundary structures. The natural action of the constructed Poisson-Lie diffeomorphism groups induces  large classes of compatible Poisson structures on ${\Bbb R}^n$, thus making it a Poisson homogeneous space. Moreover, the left-right  action of the Poisson-Lie groups $FDiff({\Bbb R}^m)\times FDiff({\Bbb R}^n)$ induces classes of compatible Poisson structures on the space $J^{\infty}({\Bbb R}^m,{\Bbb R}^n)$ of infinite jets  of smooth maps ${\Bbb R}^m\to {\Bbb R}^n$, which makes it also a Poisson homogeneous space for this action. Initial steps towards classification of these structures are taken.
\end{abstract}

\section{Introduction}

The subject of this paper is the group $FDiff({\Bbb R}^n)$ of formal diffeomorphisms  of ${\Bbb R}^n$  which we equip with various Poisson structures compatible with its group structure. The group $FDiff({\Bbb R}^n)$ is the group of $\infty$-jets of diffeomorphisms of ${\Bbb R}^n$ and as such is a subset of the set $J^{\infty}({\Bbb R}^n,{\Bbb R}^n)$ of infinite jets of smooth maps  ${\Bbb R}^n\to {\Bbb R}^n$. Lie groups with the additional structure of a Poisson manifold came to be known as Poisson-Lie groups \cite{Drinfeld:paper1}. In the infinite-dimensional setting of the group of diffeomorphisms the existence of such Poisson structures is not immediately obvious. Surprisingly, there exists a rich class of such structures. Moreover one can explicitly describe them.   We view $FDiff({\Bbb R}^n)$ as an (infinite-dimensional) formal group with a structure of an infinite-dimensional formal graded manifold.  We   consider also  its subgroup, $FDiff_0({\Bbb R}^n)\subset J^{\infty}_0({\Bbb R}^n,{\Bbb R}^n)$, of jets of smooth maps ${\Bbb R}^n\to {\Bbb R}^n$ that leave the origin $0\in{\Bbb R}^n$ fixed. The latter can be endowed with a smooth structure, for it is a  projective limit of finite-dimensional Lie groups of $k$-jets of diffeomorphisms of finite order, $k=0,1,\ldots$, that leave the origin $0\in{\Bbb R}^n$ fixed.

The motivation for studying these quasi-classical objects lies within the program of
deformation quantization \cite{BFSLF:paper, Drinfeld:IMC_report, Kontsevich:DeformationQuantization_paper, Etingof-Kazhdan:QuantizationOfLieBialgebras,Etingof-Schedler-Schiffmann:Explicit quantization}. Namely, what is sought further is the   construction of the corresponding quantum objects (quantum diffeomrphism groups, their quantum homogeneous spaces) by deforming (the structure) of their quasi-classical limits. 

An initial step towards the realization of this program was taken in \cite{Stoyanov:paper_JMP} where we formulated and began the study of the problem of
classification of Poisson-Lie structures on the infinite-dimensional
group $FDiff_0({\Bbb R}^n)$, and the 
(formal) group $FDiff({\Bbb R}^n)$ of all diffeomorphisms of ${\Bbb
R}^n$ as well as their quantization. The problem of classification of Poisson structures for
these infinite-dimensional groups as well as the classification of
their various Poisson homogeneous spaces might seem hopeless at first. However, it was completely solved in \cite{Stoyanov:paper_JMP, KS1:paper} for the groups
$FDiff_0({\Bbb R}^1)$ and $FDiff({\Bbb R}^1)$, where previously only partial results related to the existence of Lie bialgebra structures \cite{Drinfeld:paper1} on the corresponding Lie algebra  $W_1$ (Witt algebra) were known \cite{Mi:paper, Ta:paper, Dz:paper}. The moduli spaces of Poisson-Lie structures on $FDiff_0({\Bbb R}^1)$ and $FDiff({\Bbb R}^1)$ turned out to be isomorphic to ${\Bbb Z}_{+}$ and ${\Bbb Z}_{+}\cup\{-1\}$ respectively. Thus, there are
countably many equivalence classes of Poisson structures on the group
of formal diffeomorphisms of the line ${\Bbb R}^1$. 

In the present
paper we extend results obtained in \cite{Stoyanov:paper_JMP} to the general case of $FDiff_0({\Bbb R}^n)$ and $FDiff({\Bbb R}^n)$. Central to the paper are Theorems
3.1 and 3.2, which describe explicitly a large class of coboundary Poisson-Lie
structures on $FDiff({\Bbb R}^n)$ and $FDiff_0({\Bbb R}^n)$, which we also conjecture to be
exhaustive. In essence we construct solutions (skew-symmetric triangular $r$-matrices) of 
 of the Classical Yang-Baxter Equation
(CYBE) for the Lie algebra $W_n$ of formal vector fields on ${\Bbb R}^n$. We give some explicit examples and take 
initial steps towards classification of these classical $r$-matrices
and thus the corresponding Poisson-Lie groups. More precisely, we classify the $FDiff_0({\Bbb R}^n)$-orbits into which the constructed class of solutions of the CYBE splits.

The group of diffeomorphisms acts on various spaces. Two examples are the
infinite-dimensional manifold $X^{coor}$ of formal coordinate systems
in an neighborhood of a point and  the space $V_{\lambda}$ of generalized
forms $f(u^1,\ldots,u^n)(du^1\wedge\ldots\wedge du^n)^{\lambda}$. The construction of Poisson diffeomorphism groups naturally leads to the study of their Poisson actions on these spaces. Namely, one
of our objectives is to find Poisson manifolds (possibly
infinite-dimensional) which are Poisson homogeneous spaces of the Poisson-Lie  groups
$FDiff_0({\Bbb R}^n)$ and $FDiff({\Bbb R}^n)$.  For instance, in \cite{Stoyanov:paper_JMP} it was shown  that  the $FDiff_0({\Bbb R}^1)$-module $f(u)(du)^{\lambda}$
carries a large class of Poisson structures compatible with the natural action of the Poisson group $FDiff_0({\Bbb R}^1)$. In the present paper we construct two more examples of Poisson homogeneous spaces: the coordinate space ${\Bbb R}^n$ and the space of infinite jets $J^{\infty}({\Bbb R}^m,{\Bbb R}^n)$ of smooth maps ${\Bbb R}^m\to{\Bbb R}^n$.  They arise naturally from the following commutative diagram:
$$
\begin{CD} 
{\Bbb R}^m @>F>> {\Bbb R}^n \\ 
@V{X}VV @VV{Y}V\\
{\Bbb R}^m @>\overline{F}>> {\Bbb R}^n.
\end{CD}
$$
Here $X\in Diff({\Bbb R}^m)$, $Y\in Diff({\Bbb R}^n)$, $F,\overline{F}\in {\mathcal C}^{\infty}({\Bbb R}^m,{\Bbb
R}^n)$, the space of smooth maps ${\Bbb R}^m\to{\Bbb R}^n$, and where $\overline{F}=Y\circ F\circ X^{-1}$. The last defines
the left-right action of the group $Diff({\Bbb R}^n)\times Diff({\Bbb R}^m)$ on
the space ${\mathcal C}^{\infty}({\Bbb R}^m,{\Bbb R}^n)$. From these we pass to the spaces of jets of smooth maps related to the same commutative diagram: the groups of formal diffeomorphisms  $FDiff({\Bbb R}^m)$ and $FDiff({\Bbb R}^n)$, the space  of infinite jets $J^{\infty}({\Bbb R}^m,{\Bbb R}^n)$ of smooth maps ${\Bbb R}^m\to{\Bbb R}^n$, and the left-right action of
$FDiff({\Bbb R}^n)\times FDiff({\Bbb R}^m)$ on the latter space. Then we construct Poisson structures on ${\Bbb R}^m$, ${\Bbb R}^n$, and $J^{\infty}({\Bbb R}^m,{\Bbb R}^n)$ compatible with the corresponding group actions.

The paper is organized as follows.
First we construct a class of Poisson structures $\Omega_n$ on $FDiff({\Bbb R}^n)$ compatible with 
its group structure. Then we consider the action of the Poisson groups
($FDiff({\Bbb R}^n)$,$\Omega_n$)
$$
FDiff({\Bbb R}^n)\times {\Bbb R}^n\to {\Bbb R}^n, \quad u\mapsto X(u),
$$
on ${\Bbb R}^n$. We describe the class of Poisson structures on ${\Bbb R}^n$ induced by each
($FDiff({\Bbb R}^n)$,$\Omega_n$) and such that the above action
becomes Poisson. Here and in the sequel products of Poisson spaces are 
taken always with their natural product Poisson structure. Next, we
consider the action 
$$
FDiff({\Bbb R}^n)\times FDiff({\Bbb R}^m)\times J^{\infty}({\Bbb
R}^m,{\Bbb R}^n)\to J^{\infty}({\Bbb R}^m,{\Bbb R}^n), \quad F\mapsto
Y\circ F\circ X^{-1},
$$
of the Poisson groups ($FDiff({\Bbb R}^n)\times FDiff({\Bbb
R}^m),\Omega_n\times\Omega_m$) on the space of jets $J^{\infty}({\Bbb
R}^m,{\Bbb R}^n)$. It turns out
that for every pair of Poisson structures $\Omega_n$ and $\Omega_m$ there exists a
Poisson structure $\Pi$ on the infinite-dimensional space $J^{\infty}({\Bbb
R}^m,{\Bbb R}^n)$, such that
the above action is Poisson. We give an explicit description of the
corresponding Poisson structures. 

Finally, we take steps towards the solution of the classification problem. We describe the $FDiff_0({\Bbb R}^n)$-orbits (isomorphism classes)  into which the constructed set of  solutions of the CYBE (and thus the coboundary Poisson-Lie structures on $FDiff({\Bbb R}^n)$) splits. Finally we  give a short proof of the classification theorem for Poisson-Lie structures on $FDiff({\Bbb R}^1)$ \cite{Stoyanov:paper_JMP}. Some examples in low dimensions are collected in the Appendix.

{\sl Acknowledgments.} It would have been impossible to complete this work without the constant financial and moral support of Vanessa L\'opez. I am grateful to John Harnad for hospitality and
financial support at the Centre de Recherche Math\'ematique at the Universit\'e de Montr\'eal, while this paper was being written.

\section{Groups of Formal Diffeomorphisms}

Let $G_{0n}=FDiff_0({\Bbb R}^n)$ be the group of formal diffeomorphisms of  ${\Bbb R}^n$ fixing the point $0$.  As a topological group $G_{0n}$ is the group $Aut({\Bbb R}[[u^1,\ldots,u^n]])$ of all automorphisms of the  ${\Bbb R}$-algebra ${\Bbb R}[[u^1,\ldots,u^n]]$. The group $G_{0n}$ acts naturally on the space $X^{coor}$ of formal coordinate systems at $0$. The action is induced by the substitution $u\mapsto X(u)$, for $u\in {\Bbb R}^n$ and $X\in G_{0n}$. We shall identify an element
of the group $G_{0n}$ with
\be
X^i(u)=\sum_{(i_1,\ldots,i_n)\in {\Bbb Z}_+^n\backslash{\Bbb O}}x^i_{i_1i_2\ldots i_n}(u^1)^{i_1}(u^2)^{i_2}\ldots(u^n)^{i_n},\qquad i=1,\ldots,n,\label{eqn1}
\ee
where $u=(u^1,\ldots,u^n)\in {\Bbb R}^n$. Let the set $\tilde{X}^{coor}=\left\{x^i_{i_1i_2\ldots i_n}\mid (i_1,i_2,\ldots, i_n)\in {\Bbb Z}_+^n\backslash{\Bbb O},i=1,\ldots,n\right\}$ be a set of coordinates in the neighborhood of the identity (coordinates on the space of $\infty$-jets at the origin) of $G_{0n}$. Functions on this group can be identified with formal power series in $(u^1,\ldots,u^n)$ whose coefficients are smooth functions  of {\it finite} subsets of elements of $\tilde{{X}}^{coor}$, that is, functions on spaces of jets of finite order of diffeomorphisms of ${\Bbb R}^n$ fixing $0$. 

More generally, if one considers the
full group $G_n=FDiff({\Bbb R}^n)$, one can treat it as  a genuinely {\em formal} infinite-dimensional group. We identify an element of $G_n$ with
\be
X^i(u)=\sum_{(i_1,\ldots,i_n)\in {\Bbb Z}_+^n}x^i_{i_1i_2\ldots i_n}(u^1)^{i_1}(u^2)^{i_2}\ldots(u^n)^{i_n},\qquad i=1,\ldots,n.\label{eqn1}
\ee
``Functions'' on this group are formal power series in $(u^1,\ldots,u^n)$ whose coefficients are themselves formal power series in the elements  of the set of group coordinates ${X^{coor}}=\left\{x^i_{i_1i_2\ldots i_n}\mid (i_1,i_2,\ldots, i_n)\in {\Bbb Z}_+^n,i=1,\ldots,n\right\}$. More precisely,  $G_n$ (and $G_{0n}$) can be thought of as a formal graded manifold the function algebra of which is defined as follows. Assign a degree $deg$ : ${X^{coor}}\to{\Bbb Z}_+$ to each element $x^i_{i_1i_2\ldots i_n}$ of ${X^{coor}}$  by $deg(x^i_{i_1i_2\ldots i_n})=i_1+\cdots +i_n$. Consider the commutative associative algebra ${\mathcal A}={\Bbb R}[[{X^{coor}}]]$ with a unit generated by ${X^{coor}}$. The algebra ${\mathcal A}$ is naturally graded ${\mathcal A}=\oplus_{k=0}^{\infty}{\mathcal A}_k$, where  the graded components ${\mathcal A}_k$ are generated by monomials $a\in{\mathcal A}$ of degree $deg(a)=k$, for the function $deg$ can be  extended inductively to arbitrary monomials by $deg(ab)=deg(a)+deg(b)$. The unit 1 is chosen to be in ${\mathcal A}_1$. Then the algebra of functions on $G_n$ can be identified with the  graded ${\Bbb R}$-algebra ${\mathcal A}={\Bbb R}[[{X^{coor}}]]=\oplus_{k=0}^{\infty}{\mathcal A}_k$, where ${\mathcal A}_k=span\left\{a\mid deg(a)=k\right\}$. Similarly, the algebra of functions on $G_{0n}$ (if we think of it as a formal group) can be identified with $\tilde{{\mathcal A}}={\Bbb R}[[\tilde{{X}}^{coor}]]=\oplus_{k=1}^{\infty}\tilde{{\mathcal A}}_k\subset{\mathcal A}$.

As already mentioned, the multiplication in $G_n$  is defined by substitution of formal power series. It induces a comultiplication ${{\mathcal A}}\to {{\mathcal A}}\otimes{{\mathcal A}}$.  Let $X^i(u)=\sum_{i_1,\ldots,i_n}x^i_{i_1\ldots i_n}(u^1)^{i_1}\ldots(u^n)^{i_n}$ and $Y^i(u)=\sum_{i_1,\ldots,i_n}y^i_{i_1\ldots i_n}(u^1)^{i_1}\ldots(u^n)^{i_n}$ be two elements of the group $G_{n}$. The product $Z(u)$ of these two elements is defined to  be $Z(u)=X(Y(u))$:
\bea
Z^i(u)&=&\sum_{i_1,\ldots,i_n}x^i_{i_1\ldots i_n}(Y^1(u))^{i_1}\ldots(Y^n(u))^{i_n},\nonumber\\
      &=&\sum_{i_1,\ldots,i_n}z^i_{i_1\ldots i_n}(u^1)^{i_1}\ldots(u^n)^{i_n},\qquad i=1,\ldots,n,\label{eqn2}
\eea
where $z^i_{i_1\ldots i_n}=z^i_{i_1\ldots i_n}(x,y)\in{{\mathcal A}}\otimes{{\mathcal A}}$ (here we have identified $x^i_{i_1\ldots i_n}$ with $x^i_{i_1\ldots i_n}\otimes 1$ and $y^i_{i_1\ldots i_n}$ with $1\otimes x^i_{i_1\ldots i_n}$). These are polynomials in  a finite number of $x^i_{i_1\ldots i_n}$'s and $y^i_{i_1\ldots i_n}$'s in the case of $G_{0n}$. In the case of $G_n$, $z^i_{i_1\ldots i_n}(x,y)\in\oplus_{k=0}^{\infty}{\mathcal A}_k\otimes{\mathcal A}_{i_1+\ldots+i_n}\subset{\mathcal A}\otimes{\mathcal A}$ are formal power series (defining a formal group multiplication law \cite{Se:paper}) in $x^i_{i_1\ldots i_n}$ and $y^i_{i_1\ldots i_n}$.  Namely,
\be
z^i_{i_1\ldots i_n}=\sum_{k_1=0}^{\infty}\ldots\sum_{k_n=0}^{\infty}x^i_{k_1\ldots k_n}\sum_{(\stackrel{\mbox{\tiny all partitions}}{s_{11}+\cdots+s_{1n}=i_1})}\ldots\sum_{(\stackrel{\mbox{\tiny all partitions}}{s_{n1}+\cdots+s_{nn}=i_n})}Y_{s_{11}\ldots s_{1n}}^{(k_1)}\ldots Y_{s_{n1}\ldots s_{nn}}^{(k_n)},
\ee
where
\be
Y_{s_1\ldots s_n}^{(k_p)}=\sum_{(\stackrel{\mbox{\tiny all partitions}}{\alpha_{11}+\cdots+\alpha_{1k_p}=s_1})}\ldots\sum_{(\stackrel{\mbox{\tiny all partitions}}{\alpha_{n1}+\cdots+\alpha_{nk_p}=s_n})}y_{\alpha_{11}\ldots\alpha_{n1}}^{p}\ldots y_{\alpha_{1k_p}\ldots\alpha_{nk_p}}^{p},
\ee
and where $deg\left(Y_{s_1\ldots s_n}^{(k_p)}\right)=s_1+\cdots +s_n$. The summation $\sum_{(\stackrel{\mbox{\tiny all partitions}}{\alpha_{1}+\cdots+\alpha_{k}=s})}$ is over all partitions of $s\in{\Bbb Z}_{+}$ into a sum of $k$ non-negative integers.

The inverse $X\mapsto X^{-1}$ of an element $X$ of $G_n$ is defined by $X^i(X^{-1}(u))=u^i=(X^{-1})^i(X(u))$. It defines an antipode ${{\mathcal A}}\to {{\mathcal A}}$ on the algebra ${{\mathcal A}}$. For the existence of the inverse it is necessary and sufficient that the $n\times n$ matrix
\be
X_0=\pmatrix{x^1_{10\ldots 0}&x^1_{01\ldots 0}&\ldots&x^1_{00\ldots 1}\cr x^2_{10\ldots 0}&x^2_{01\ldots 0}&\ldots&x^2_{00\ldots 1}\cr \vdots&\vdots&\ddots&\vdots\cr x^n_{10\ldots 0}&x^n_{01\ldots 0}&\ldots&x^n_{00\ldots 1}}\in \mbox{Mat}(n,{{\mathcal A}}),\label{eqn3V3}
\ee
is invertible, that is, its determinant $\mbox{det}X_0$ is an invertible element in ${{\mathcal A}}$. We make this a defining condition for the group $G_n$ (and $G_{0n}$) and assume that $(\mbox{det}X_0)^{-1}\in{{\mathcal A}_0} $.
\begin{example} Consider  the simplest case of $FDiff_0({\Bbb R}^1)$ and $FDiff({\Bbb R}^1)$. It is not difficult to see that the formulae for the multiplication  are given by
\be
z_i(x,y) =\sum_{k=1}^ix_k\sum_{(\sum_{\alpha=1}^{k}{s_{\alpha}})=i}{y_{s_1}\ldots y_{s_k}}\quad \in \oplus_{k=1}^{i}\tilde{{\mathcal A}}_k\otimes{\tilde{{\mathcal A}}_i}\subset{\tilde{{\mathcal A}}}\otimes{\tilde{{\mathcal A}}}, \quad i\ge 1,
\ee
for $FDiff_0({\Bbb R}^1)$ and 
\begin{eqnarray}z_0(x,y)&=&\sum_{k=0}^{\infty}x_ky_0^k\quad \in \oplus_{k=0}^{\infty}{{\mathcal A}}_k\otimes{{{\mathcal A}}_0}\subset{{{\mathcal A}}}\otimes{{{\mathcal A}}},\label{eqnVII2}\\
z_i(x,y)&=&\sum_{k=1}^{\infty}x_k\sum_{\bigl(\sum_{\alpha=1}^{k}s_{\alpha}\bigr)=i}y_{s_1}\ldots y_{s_{k}}\quad \in \oplus_{k=1}^{\infty}{{\mathcal A}}_k\otimes{{{\mathcal A}}_i}\subset{{{\mathcal A}}}\otimes{{{\mathcal A}}}, \quad   i\ge 1.\label{eqnVII3}
\end{eqnarray}
for $FDiff({\Bbb R}^1)$ respectively.
The (recursively computed) first several formulae that determine the inverse $X^{-1}(u)=\sum_{k=1}^{\infty}\overline{x}_ku^k$ of an element $X\in FDiff_0({\Bbb R}^1)$ are 
\bea
\overline{x}_1&=&\frac{1}{x_1},\nonumber\\
\overline{x}_2&=&-\frac{x_2}{x_1^3},\nonumber\\
\overline{x}_3&=&-\frac{x_3}{x_1^4}+\frac{2x_2^2}{x_1^5},\nonumber\\
&\vdots&\nonumber
\eea
and those that determine the inverse $X^{-1}(u)=\sum_{k=0}^{\infty}\overline{x}_ku^k$ of an element $X\in FDiff({\Bbb R}^1)$ are given by
\bea
\overline{x}_0&=&-\frac{x_0}{x_1}-\frac{x_0^2x_2}{x_1^3}+\frac{x_0^3x_3}{x_1^4}-\frac{1}{x_1^5}(2x_0^3x_2^2+x_0^4x_4)+{\bf O}(x_1^{-6}),\nonumber\\
\overline{x}_1&=&\frac{1}{x_1}+\frac{2x_0x_2}{x_1^3}-\frac{3x_0^2x_3}{x_1^4}+\frac{1}{x_1^5}(6x_0^2x_2^2+4x_0^3x_4)+{\bf O}(x_1^{-6}),\nonumber\\
\overline{x}_2&=&-\frac{x_2}{x_1^3}+\frac{3x_0x_3}{x_1^4}-\frac{1}{x_1^5}(6x_0x_2^2+6x_0^2x_4)+{\bf O}(x_1^{-6}),\nonumber\\
\overline{x}_3&=&-\frac{x_3}{x_1^4}+\frac{1}{x_1^5}(2x_2^2+4x_0x_4)+{\bf O}(x_1^{-6}),\nonumber\\
&\vdots&\nonumber
\eea
where ${\bf O}(x_1^{-6})\in x^{-6}_1{\Bbb R}[[{X^{coor}}\backslash\{x_1\},x^{-1}_1]]$.
\end{example}
We conclude this section with the derivation of  several useful formulae stemming immediately from the above considerations. We do this for the group $G_{0n}$. The formulae for the group $G_n$ can be derived analogously. The first two are
\be
\frac {\partial Z^i(u)}{\partial x^j_{j_1\ldots j_n}}=\delta^i_j(Y^1(u))^{j_1}\ldots(Y^n(u))^{j_n},\label{eqn3}
\ee
and
\be
\frac {\partial Y^i(u)}{\partial y^j_{j_1\ldots j_n}}=\delta^i_j(u^1)^{j_1}\ldots(u^n)^{j_n}.\label{eqn4}
\ee
Then, we also have
\bea
\frac {\partial Z^i(u)}{\partial y^j_{j_1\ldots j_n}}&=&\sum_{i_1,\ldots,i_n}\sum_{k=1}^{n}i_kx^i_{i_1\ldots i_n}\frac {\partial Y^k(u)}{\partial y^j_{j_1\ldots j_n}}(Y^1(u))^{i_1}\ldots(Y^{k}(u))^{i_{k}-1}\ldots(Y^n(u))^{i_n}\nonumber\\
&=&\sum_{i_1,\ldots,i_n}\sum_{k=1}^{n}i_kx^i_{i_1\ldots i_n}(u^1)^{j_1}\ldots(u^n)^{j_n}(Y^1(u))^{i_1}\ldots(Y^{k}(u))^{i_{k}-1}\ldots(Y^n(u))^{i_n}\nonumber\\
&=&(u^1)^{j_1}\ldots(u^n)^{j_n}\left(X_{*Y(u)}\right)^i_j,\label{eqn5}
\eea
where $X_{*Y(u)}$ is the derivative of the map $X$ at the point $Y(u)$ and $\delta^i_j$ stands for the Kronecker symbol.

\section{Poisson-Lie Structures on $FDiff_0({\Bbb R}^n)$ and $FDiff({\Bbb R}^n)$}

 Our task in this section is  to find out whether the groups
introduced in the preceding section can be equipped with Poisson structure compatible with the group structure \cite{Drinfeld:paper1}. We shall consider first the group $G_{0n}$ assuming that it is equipped with a smooth structure, the space of smooth functions $C^{\infty}(G_{0n})$ being the inductive limit of spaces of smooth functions on groups of finite jets of diffeomorphisms \cite{Stoyanov:paper_JMP}.  In other words, we seek to
construct a Poisson tensor $\omega$ with infinitely many components in $C^{\infty}(G_{0n})$. If on the contrary the group  $G_{0n}$  is interpreted as a formal group then $\omega$ can be interpreted as a bi-derivation in the algebra $\tilde{{\mathcal A}}\otimes\tilde{{\mathcal A}}$ where the components of $\omega$ belong to $\tilde{{\mathcal A}}$. This is the case of $G_n$ which we think of as a formal group the function algebra on which is the infinitely generated algebra ${\mathcal A}$ (cf. with the previous section). The tensor  $\omega$ has to be such that the multiplication map $G_{0n}\times G_{0n}\to G_{0n}$ is Poisson. Here the space $G_{0n}\times G_{0n}$  is equipped with the product Poisson structure. In order to express the latter condition we shall introduce a new object. Let ${\bf i}=(i_1,\ldots,i_n)\in{\Bbb Z}_{+}^n\backslash{\Bbb O}$ denote a multi-index. We shall write $x^{\{i,{\bf i}\}}$ for $x^i_{i_1\ldots i_n}$. 

A Poisson tensor on the group $G_{0n}$, which in general depends on the infinite jets of a group element $X$ at two points $u,v\in {\Bbb R}^n$, is defined by $n^2$ power series $\Omega^{ij}(j^{\infty}_{u}X,j^{\infty}_{v}X)\in C^{\infty}(G_{0n})[[u^1,\ldots,u^n,v^1,\ldots,v^n]]$:
\be
%\left\{X^i(u),X^j(v)\right\}&:=&\Omega^{ij}(j^{\infty}_{u}X,j^{\infty}_{v}X)\nonumber\\
\Omega^{ij}(j^{\infty}_{u}X,j^{\infty}_{v}X):=\sum_{i_1,\ldots,i_n}\sum_{j_1,\ldots,j_n}\omega^{\{i,i_1,\ldots,i_n\}\{j,j_1,\ldots,j_n\}}(x)(u^1)^{i_1}\ldots(u^n)^{i_n}(v^1)^{j_1}\ldots(v^n)^{j_n}.\label{eqn6}
\ee
 We shall use also the representation $\left\{X^i(u),X^j(v)\right\}\equiv \Omega^{ij}(j^{\infty}_{u}X,j^{\infty}_{v}X)$. The latter follows from   the fact that we  have $\left\{x^i_{i_1,\ldots,i_n},x^j_{j_1,\ldots,j_n}\right\}=\omega^{\{i,i_1,\ldots,i_n\}\{j,j_1,\ldots,j_n\}}(x)$ for the Poisson brackets between coordinates.

Thus, a Poisson structure on $G_{0n}$ is determined by a tensor
$\omega$ with infinitely many components
$\omega^{\{i,i_1,\ldots,i_n\}\{j,j_1,\ldots,j_n\}}(x)$ which are
smooth functions on the spaces of jets of a finite order of
diffeomorphisms of $\Bbb R^n$ fixing the point $0$, or in the case of $G_n$ -- formal power series which  belong to  ${\mathcal A}$. We shall write the
components of the tensor $\omega$ as $\omega^{\{i,{\bf i}\}\{j,{\bf j}\}}$, where
$i,j\in \{1,\ldots,n\}$ and ${\bf i},{\bf j}\in{\Bbb Z}_{+}^n\backslash{\Bbb O}$ run over an infinite set of
$n$-tuples of non-negative integers with the one consisting only of zeros excluded. 
\begin{definition}
A Poisson-Lie structure on $G_{0n}$ is defined by a ``tensor'' $\Omega^{ij}(j^{\infty}_{u}X,j^{\infty}_{v}X)$ having the properties:\newline
(a) $\Omega^{ij}(j^{\infty}_{u}X,j^{\infty}_{v}X)=-\Omega^{ji}(j^{\infty}_{v}X,j^{\infty}_{u}X)$;\newline
(b) $\Omega^{ij}(j^{\infty}_{u}X,j^{\infty}_{v}X)$ satisfies the $n^2$ functional equations:
\be
\Omega^{ij}(j^{\infty}_{u}Z,j^{\infty}_{v}Z)=\Omega^{ij}(j^{\infty}_{Y(u)}X,j^{\infty}_{Y(v)}X)+\sum_{k=1}^n\sum_{l=1}^n\left(X_{*Y(u)}\right)^i_k\left(X_{*Y(v)}\right)^j_l\Omega^{kl}(j^{\infty}_{u}Y,j^{\infty}_{v}Y);\label{eqn8}
\ee
(c) $\Omega^{ij}(j^{\infty}_{u}X,j^{\infty}_{v}X)\equiv \left\{X^i(u),X^j(v)\right\}$ satisfies the Jacobi identities, which we write as
\be
\{\{X^i(u),X^j(v)\},X^k(w)\}+\{\{X^k(w),X^i(u)\},X^j(v)\}+\{\{X^j(v),X^k(w)\},X^i(u)\}=0,\label{eqn9}
\ee
and which are to be understood as $n^3$ differential identities holding
in $C^{\infty}(G_{0n})[[u,v,w]]$. The last two terms in (\ref{eqn9})
are obtained by cyclicly permuting the pairs $(i,u),(j,v),(k,w)$ in
the first. In the above formulae $u,v,w$ stand for $u=(u^1,\ldots,u^n)$,
$v=(v^1,\ldots,v^n)$, and $w=(w^1,\ldots,w^n)$.

Similarly for the group $G_n$, in which case (\ref{eqn9}) holds in ${\mathcal A}[[u,v,w]]$.
\end{definition}

The condition (b) in the definition of a Poisson-Lie group is, in fact, the requirement that the multiplication map is Poisson. To see this we note that the latter  is equivalent to the infinite system of functional
equations
\be
\omega^{\{i,{\bf i}\}\{j,{\bf j}\}}(z)=\sum_{{\bf k},{\bf l}}\sum_{k,l=1}^n\left[\frac{\partial z^{\{i,{\bf i}\}}}{\partial x^{\{k,{\bf k}\}}}\frac{\partial z^{\{j,{\bf j}\}}}{\partial x^{\{l,{\bf l}\}}}\omega^{\{k,{\bf k}\}\{l,{\bf l}\}}(x)+\frac{\partial z^{\{i,{\bf i}\}}}{\partial y^{\{k,{\bf k}\}}}\frac{\partial z^{\{j,{\bf j}\}}}{\partial y^{\{l,{\bf l}\}}}\omega^{\{k,{\bf k}\}\{l,{\bf l}\}}(y)\right]\label{eqn7}
\ee
for the components of the tensor $\omega$ \cite{Stoyanov:paper_JMP}.
Note that the sums over ${\bf k},{\bf l}$ above are {\it finite}, since for $G_{0n}$ the $z$'s all depend on a finite number of the $x$'s and $y$'s. For $G_n$ this is not the case and equations (\ref{eqn7}) should be interpreted as holding in the graded algebra ${\mathcal A}$. Now we recast this infinite system of functional equations in the language of power series. This is the content of the following proposition.
\begin{proposition}
The multiplicativity conditions (\ref{eqn7}) are equivalent to the $n^2$ functional equations
\be
\Omega^{ij}(j^{\infty}_{u}Z,j^{\infty}_{v}Z)=\Omega^{ij}(j^{\infty}_{Y(u)}X,j^{\infty}_{Y(v)}X)+\sum_{k=1}^n\sum_{l=1}^n\left(X_{*Y(u)}\right)^i_k\left(X_{*Y(v)}\right)^j_l\Omega^{kl}(j^{\infty}_{u}Y,j^{\infty}_{v}Y).\label{eqn8}
\ee
\end{proposition}
{\it Proof:} Let us multiply both sides of (\ref{eqn7}) by the monomial $(u^1)^{i_1}\ldots(u^n)^{i_n}(v^1)^{j_1}\ldots(v^n)^{j_n}$ and sum over the sets of $n$-tuples $(i_1,\ldots,i_n)$ and $(j_1,\ldots,j_n)$. By the definition of
$\Omega$ the left hand side of (\ref{eqn7}) becomes
\bdm
\sum_{i_1,\ldots,i_n}\sum_{j_1,\ldots,j_n}\omega^{\{i,i_1,\ldots,i_n\}\{j,j_1,\ldots,j_n\}}(z)(u^1)^{i_1}\ldots(u^n)^{i_n}(v^1)^{j_1}\ldots(v^n)^{j_n}=\Omega^{ij}(j^{\infty}_{u}Z,j^{\infty}_{v}Z).
\edm
Thus, after using formulae (\ref{eqn3}) and (\ref{eqn4}) we obtain
\bea
\Omega^{ij}(j^{\infty}_{u}Z,j^{\infty}_{v}Z)&=&\sum_{{\bf k},{\bf l}}\sum_{k,l=1}^n\left[\frac{\partial Z^i(u)}{\partial x^{\{k,{\bf k}\}}}\frac{\partial Z^j(v)}{\partial x^{\{l,{\bf l}\}}}\omega^{\{k,{\bf k}\}\{l,{\bf l}\}}(x)+\frac{\partial Z^i(u)}{\partial y^{\{k,{\bf k}\}}}\frac{\partial Z^j(v)}{\partial y^{\{l,{\bf l}\}}}\omega^{\{k,{\bf k}\}\{l,{\bf l}\}}(y)\right]\nonumber\\
&=&\Omega^{ij}(j^{\infty}_{Y(u)}X,j^{\infty}_{Y(v)}X)+\sum_{k=1}^n\sum_{l=1}^n\left(X_{*Y(u)}\right)^i_k\left(X_{*Y(v)}\right)^j_l\Omega^{kl}(j^{\infty}_{u}Y,j^{\infty}_{v}Y).\nonumber
\eea
This completes the proof. \hfill \fbox{}

In what follows we shall assume that summation over repeated indices is understood and drop the summation signs unless it is explicitly
stated otherwise. 

A class of solutions of the system of functional equations (\ref{eqn8}) is described by the following proposition.
\begin{proposition} For any set of $n^2$ formal power series $\varphi^{ij}(u,v)\in{\Bbb R}[[u,v]]$ with the property $\varphi^{ij}(u,v)=-\varphi^{ji}(v,u)$ the functional equation (\ref{eqn8}) admits 
\be
\Omega^{ij}(j^{\infty}_{u}X,j^{\infty}_{v}X)=\left(X_{*u}\right)^i_k\left(X_{*v}\right)^j_l\varphi^{kl}(u,v)-\varphi^{ij}(X(u),X(v)),\label{eqn10}
\ee
as solution.
\end{proposition}
\pf\ The proof consists of a direct verification. 
We substitute the expression for $\Omega$ in both sides of (\ref{eqn8}). Then the left hand side reads
$$
\left(Z_{*u}\right)^i_k\left(Z_{*v}\right)^j_l\varphi^{kl}(u,v)-\varphi^{ij}(Z(u),Z(v))\stackrel{\rm ?}{=}
$$
and the right hand side becomes
\bea
\lefteqn{\stackrel{\rm ?}{=}\left(X_{*Y(u)}\right)^i_k\left(X_{*Y(v)}\right)^j_l\varphi^{kl}(Y(u),Y(v))-\varphi^{ij}(X(Y(u)),X(Y(v)))+}\nonumber\\
& &\mbox{\phantom{$\Phi$}} +\left(X_{*Y(u)}\right)^i_k\left(X_{*Y(v)}\right)^j_l\left(Y_{*u}\right)^k_s\left(Y_{*v}\right)^l_p\varphi^{sp}(u,v)-\nonumber\\
& &\mbox{\phantom{$\Phi$}}-\left(X_{*Y(u)}\right)^i_k\left(X_{*Y(v)}\right)^j_l\varphi^{kl}(Y(u),Y(v)).\nonumber
\eea
Since $Z=X\circ Y$, we have $Z_{*u}=X_{*Y(u)}Y_{*u}$, and  we see
immediately that we obtain an identity. \epf

\begin{remark} Formula (\ref{eqn10}) can be interpreted as follows. Let $W_n$ denote the algebra of formal vector fields on ${\Bbb R^n}$.  This algebra has an interpretation as the algebra of $\infty$-jets of smooth vector fields at the origin equipped with the projective limit topology.  Consider the space of maps ${\mathcal M}=Map({\Bbb R^n}\times{\Bbb R^n},W_n\widehat{\otimes}W_n)$, where $W_n\widehat{\otimes}W_n$ is the completed tensor product. This space is a $G_{0n}$-module. The action is given by the right hand side of formula (\ref{eqn10}). Thus the solution (\ref{eqn10}) is a 1-coboundary in $H^1(G_{0n},{\mathcal M})$.
Later we shall see that $\varphi^{ij}(u,v)$ is nothing but a skew-symmetric ($\varphi^{ij}(u,v)=-\varphi^{ji}(v,u)$) classical $r-matrix$ \cite{STS2:paper} for the Lie algebra $W_n$.
\end{remark}

The Jacobi identities (\ref{eqn9}) imply what set of equations the $n^2$ series $\varphi^{ij}(u,v)$ must satisfy in order for $\Omega^{ij}(j^{\infty}_{u}X,j^{\infty}_{v}X)$ to define a Poisson structure. These are described by the proposition that follows.
\begin{proposition} The tensor $\Omega(j^{\infty}_{u}X,j^{\infty}_{v}X)$ given by equation (\ref{eqn10}) defines a Poisson structure  on $G_{0n}$ compatible with the group structure if and only if the following set of $n^3$ equations hold:
\be
\Phi^{ijk}(X(u),X(v),X(w))=\left(X_{*u}\right)^i_s\left(X_{*v}\right)^j_q\left(X_{*w}\right)^k_t\Phi^{sqt}(u,v,w),\label{eqn11}
\ee
where $i,j,k=1,\ldots,n$, and 
\bea
\lefteqn{\Phi^{ijk}(u,v,w)=\varphi^{ks}(w,u)\frac{\partial \varphi^{ij}(u,v)}{\partial u^s}+\varphi^{sk}(v,w)\frac{\partial \varphi^{ji}(v,u)}{\partial v^s}+}\nonumber\\
& & \mbox{\phantom{$\Phi^{ijk}(u,v,w)=$}} + \varphi^{is}(u,v)\frac{\partial \varphi^{jk}(v,w)}{\partial v^s}+\varphi^{si}(w,u)\frac{\partial \varphi^{kj}(w,v)}{\partial w^s}+\nonumber\\
& & \mbox{\phantom{$\Phi^{ijk}(u,v,w)=$}} + \varphi^{js}(v,w)\frac{\partial \varphi^{ki}(w,u)}{\partial w^s}+\varphi^{sj}(u,v)\frac{\partial \varphi^{ik}(u,w)}{\partial u^s}.\label{eqn12}
\eea
In other words, $\Omega(j^{\infty}_{u}X,j^{\infty}_{v}X)$ satisfies the Jacobi identities (\ref{eqn9}) if and only if $\Phi(u,v,w)$ is an invariant for the action of the group $G_{0n}$, $\Phi\in \mbox{Inv}_{G_{0n}}(W_n\widehat{\otimes}W_n\widehat{\otimes}W_n)$.
\end{proposition}
\pf We give a sketch of the calculations involved. Namely, we compute the first term in (\ref{eqn9}) with the use of formula (\ref{eqn10}):
$$
\{\{X^i(u),X^j(v)\},X^k(w)\}=\varphi^{sp}(u,v)\left\{\frac{\partial X^i(u)}{\partial u^s}\frac{\partial X^j(v)}{\partial v^p},X^k(w)\right\}-\left\{\varphi^{ij}(X(u),X(v)),X^k(w)\right\}
$$
$$
=\varphi^{sp}(u,v)\frac{\partial X^i(u)}{\partial u^s}\frac{\partial }{\partial v^p}\left\{X^j(v),X^k(w)\right\}+\varphi^{sp}(u,v)\frac{\partial X^j(v)}{\partial v^p}\frac{\partial }{\partial u^s}\left\{X^i(u),X^k(w)\right\}
$$
$$
-\partial_{1s}\varphi^{ij}(X(u),X(v))\left\{X^s(u),X^k(w)\right\}-\partial_{2s}\varphi^{ij}(X(u),X(v))\left\{X^s(v),X^k(w)\right\}
$$
$$
=\varphi^{sp}(u,v)\frac{\partial X^i(u)}{\partial u^s}\left[\frac{\partial^2X^j(v)}{\partial v^p\partial v^q}\frac{\partial X^k(w)}{\partial w^t}\varphi^{qt}(v,w)+\frac{\partial X^j(v)}{\partial v^q}\frac{\partial X^k(w)}{\partial w^t}\partial_{1p}\varphi^{qt}(v,w)\right]+
$$
$$
+\varphi^{sp}(u,v)\frac{\partial X^j(v)}{\partial v^p}\left[\frac{\partial^2X^i(u)}{\partial u^s\partial u^q}\frac{\partial X^k(w)}{\partial w^t}\varphi^{qt}(v,w)+\frac{\partial X^i(u)}{\partial u^q}\frac{\partial X^k(w)}{\partial w^t}\partial_{1s}\varphi^{qt}(u,w)\right]+
$$
$$
-\varphi^{sp}(u,v)\frac{\partial X^i(u)}{\partial u^s}\frac{\partial X^q(v)}{\partial v^p}\partial_{1q}\varphi^{jk}(X(v),X(w))-\varphi^{sp}(u,v)\frac{\partial X^j(v)}{\partial v^p}\frac{\partial X^q(u)}{\partial u^s}\partial_{1q}\varphi^{ik}(X(u),X(w))
$$
$$
-\partial_{1s}\varphi^{ij}(X(u),X(v))\frac{\partial X^s(u)}{\partial u^q}\frac{\partial X^k(w)}{\partial w^t}\varphi^{qt}(u,w)+\partial_{1s}\varphi^{ij}(X(u),X(v))\varphi^{sk}(X(u),X(w))
$$
$$
-\partial_{2s}\varphi^{ij}(X(u),X(v))\frac{\partial X^s(v)}{\partial v^q}\frac{\partial X^k(w)}{\partial w^t}\varphi^{qt}(u,w)+\partial_{2s}\varphi^{ij}(X(u),X(v))\varphi^{sk}(X(v),X(w)).
$$
Here $\partial_{1s}$ and $\partial_{2s}$ denote derivatives with respect to the first and the second arguments. The formulae for the last two terms in (\ref{eqn9}) can be obtained from the formula above by a cyclic
permutation of the pairs $(i,u)$, $(j,v)$ and $(k,w)$. Adding the resulting expressions we obtain (\ref{eqn11}) after cancellations. \hfill \fbox{}

In particular, the equations (\ref{eqn11}) hold whenever $\Phi^{ijk}(u,v,w)=0$. These comprise the Classical Yang-Baxter Equation for the Lie algebra of formal vector fields $W_n$. Therefore before turning the discussion to solutions of these equations we shall introduce the corresponding algebraic notions for the benefit of a reader more familiar with
the language of Lie bialgebras. 

Let ${\mathcal X}^i(u)\frac{\partial}{\partial u^i}$ and ${\mathcal Y}^i(u)\frac{\partial}{\partial u^i}$ be two formal vector fields in $W_n$. Let their Lie bracket be $[{\mathcal X},{\mathcal Y}]^i(u)\frac{\partial}{\partial u^i}$ where 
\be
[{\mathcal X},{\mathcal Y}]^i(u)={\mathcal X}^i(u)\frac{\partial {\mathcal Y}^j(u)}{\partial u^i}-{\mathcal Y}^i(u)\frac{\partial {\mathcal X}^j(u)}{\partial u^i}.\label{eqn13}
\ee
Let $\delta : W_n\to W_n\widehat{\otimes} W_n$ be a coalgebra map, from $W_n$ to the completed tensor product $W_n\widehat{\otimes} W_n$.
It maps a vector field to a bivector field 
\be
{\mathcal X}^i(u)\frac{\partial}{\partial u^i}\stackrel{\delta}{\longmapsto} \varphi^{ij}_{\mathcal X}(u,v)\frac{\partial}{\partial u^i}\widehat{\otimes}\frac{\partial}{\partial v^j}, 
\ee
where $\varphi^{ij}_{\mathcal X}(u,v)\in {\Bbb R}[[u^1,\ldots,u^n,v^1,\ldots,v^n]]$ are formal power series which generally depend on
$(u,{\mathcal X}(u),{\mathcal X}'(u),\ldots;v,{\mathcal X}(v),{\mathcal X}'(v),\ldots)$, i.e., the infinite jets of ${\mathcal X}$ at the points $u$ and $v$.

We recall that the map $\delta$ defines a Lie bialgebra structure on $W_n$ if and only if the following three conditions are satisfied:
\begin{eqnarray}(a)\qquad && \tau\circ\delta=-\delta\nonumber\\
(b) \qquad && \delta\bigl([{\mathcal X},{\mathcal Y}]\bigr)=ad_{\mathcal X}\delta({\mathcal Y})-ad_{\mathcal Y}\delta({\mathcal X}),\ \ \ {\mathcal X},{\mathcal Y}\in W_n,\nonumber\\
(c)\qquad && [1\widehat{\otimes} 1\widehat{\otimes} 1+(\tau\widehat{\otimes}  1)(1\widehat{\otimes}\tau)+(1\widehat{\otimes}\tau)(\tau\widehat{\otimes} 1)](1\widehat{\otimes}\delta)\circ\delta=0 ,\nonumber
\end{eqnarray}
where $\tau$ is the transposition map $\tau{\ :\ }{W_n}\widehat{\otimes}{W_n}\to{W_n}\widehat{\otimes}{W_n}$ defined by $\tau({\mathcal X}\widehat{\otimes} {\mathcal Y})={\mathcal Y}\widehat{\otimes} {\mathcal X}$, for any ${\mathcal X},{\mathcal Y}\in{W_n}$. The second and first conditions mean that $\delta$ is a Lie algebra 1-cocycle in $H^1(W_n,W_n\widehat{\otimes} W_n)$. The third condition is the so called co-Jacobi identity and $\delta$ satisfying it means that $\delta$ is coassociative. 

The condition (a) is equivalent to $\varphi^{ij}_{\mathcal X}(u,v)=-\varphi^{ji}_{\mathcal X}(v,u)$.
Further, computing the action
\be
ad_{\mathcal X}\delta({\mathcal Y})=\left[{\mathcal X}^i(u)\frac{\partial \varphi^{kl}_{{\mathcal Y}}}{\partial u^i}-\varphi^{il}_{{\mathcal Y}}(u,v)\frac{\partial {\mathcal X}^k}{\partial u^i}+{\mathcal X}^i(v)\frac{\partial \varphi^{kl}_{{\mathcal Y}}}{\partial v^i}-\varphi^{ki}_{{\mathcal Y}}(u,v)\frac{\partial {\mathcal X}^l}{\partial v^i}\right]\frac{\partial }{\partial u^k}\widehat{\otimes}\frac{\partial }{\partial u^l},\label{eqn15}
\ee
for any two vector fields ${\mathcal X},{\mathcal Y}\in W_n$, it is not difficult to see that the 1-cocycle condition for $W_n$ is equivalent to the following  system of equations which must hold in ${\Bbb R}[[u^1,\ldots,u^n,v^1,\ldots,v^n]]$:
\bea
\varphi^{kl}_{[{\mathcal X},{\mathcal Y}]}(u,v)&=&{\mathcal X}^i(u)\frac{\partial \varphi^{kl}_{{\mathcal Y}}}{\partial u^i}-\varphi^{il}_{{\mathcal Y}}(u,v)\frac{\partial {\mathcal X}^k}{\partial u^i}+{\mathcal X}^i(v)\frac{\partial \varphi^{kl}_{{\mathcal Y}}}{\partial v^i}-\varphi^{ki}_{{\mathcal Y}}(u,v)\frac{\partial {\mathcal X}^l}{\partial v^i}\nonumber\\
&&-{\mathcal Y}^i(u)\frac{\partial \varphi^{kl}_{{\mathcal X}}}{\partial u^i}+\varphi^{il}_{{\mathcal X}}(u,v)\frac{\partial {\mathcal Y}^k}{\partial u^i}-{\mathcal Y}^i(v)\frac{\partial \varphi^{kl}_{{\mathcal X}}}{\partial v^i}+\varphi^{ki}_{{\mathcal X}}(u,v)\frac{\partial {\mathcal Y}^l}{\partial v^i}.\label{eqn16}
\eea

In what follows we shall analyze the class of coboundary solutions of the above system of equations. Namely, we shall consider  the ones for which $\varphi_{\mathcal X}=ad_{\mathcal X}\varphi$, where $\varphi=\varphi^{ij}(u,v)\frac{\partial }{\partial u^i}\widehat{\otimes}\frac{\partial }{\partial u^j}$ is some formal bi-vector field with $\varphi^{ij}(u,v)\in{\Bbb R}[[u,v]]$ and such that $\varphi^{ij}(u,v)=-\varphi^{ji}(v,u)$. In other words, the $n^2$ series $\varphi^{ij}(u,v)$ are the elements of the classical skew-symmetric $r$-matrix for $W_n$. The class of coboundary solutions is given by
\be
\varphi^{ij}_{\mathcal X}(u,v)={\mathcal X}^k(u)\frac{\partial \varphi^{ij}}{\partial u^k}-\varphi^{kj}\frac{\partial {\mathcal X}^i(u)}{\partial u^k}+{\mathcal X}^k(v)\frac{\partial \varphi^{ij}}{\partial v^k}-\varphi^{ik}\frac{\partial {\mathcal X}^j(v)}{\partial v^k}.\label{eqn17}
\ee
With $\varphi^{ij}_{\mathcal X}(u,v)$ given by formula (\ref{eqn17}) it is now not difficult to establish that the requirement for coassociativity of $\delta$ is equivalent to the following system of $n^3$ linear partial differential equations 
\be
{\mathcal X}^s(u)\frac{\partial \Phi^{ijk}}{\partial u^s}+{\mathcal X}^s(v)\frac{\partial \Phi^{ijk}}{\partial v^s}+{\mathcal X}^s(w)\frac{\partial \Phi^{ijk}}{\partial w^s}-\frac{\partial {\mathcal X}^i(u)}{\partial u^s}\Phi^{sjk}-\frac{\partial {\mathcal X}^j(v)}{\partial v^s}\Phi^{isk}-\frac{\partial {\mathcal X}^k(w)}{\partial w^s}\Phi^{ijs}=0,\label{eqn18}
\ee
which must hold in ${\Bbb R}[[u,v,w]]$. Here $\Phi$ is given by
(\ref{eqn12}).  We obtain that the coassociativity of $\delta$
implies that $ad_{\mathcal X}(\Phi)=0$, i.e., $\Phi\in\mbox{Inv}_{W_n}(W_n\widehat{\otimes}W_n\widehat{\otimes}W_n)$.
Equation (\ref{eqn18}) is what is known as the generalized classical
Yang-Baxter equation (GCYBE) for the Lie algebra $W_n$ (in our
context). It is the infinitesimal part of (\ref{eqn11}) and is obtained from equation (\ref{eqn11}) by writing $X^i(u)=u^i+\epsilon{\mathcal X}^i(u)$, $i=1,\ldots,n$, for some diffeomorphism ${\mathcal X}$, expanding around the identity diffeomorphism and keeping only terms of order $\epsilon$.

In general, in order to construct Poisson structures one needs to find the space of invariants $\mbox{Inv}_{G_{0n}}(W_n\widehat{\otimes}W_n\widehat{\otimes}W_n)$ and then construct solutions of the system of functional partial differential equations (\ref{eqn12}). An important non-trivial problem is the
classification of these Poisson structures. It has been solved in a few cases, all for finite-dimensional groups (notably the classification result for complex finite dimensional semi-simple Lie algebras \cite{BD:paper} and the explicit quantization in \cite{Etingof-Schedler-Schiffmann:Explicit quantization}) with the exception of the case of $G_1$ for which the classification problem was completely solved in \cite{Stoyanov:paper_JMP} (see below). 

In the present paper we shall concentrate on solutions of equation (\ref{eqn12}) with $\Phi=0$, i.e., solutions of the classical Yang-Baxter equation 
\bea
\lefteqn{\varphi^{ks}(w,u)\frac{\partial \varphi^{ij}(u,v)}{\partial u^s}+\varphi^{sk}(v,w)\frac{\partial \varphi^{ji}(v,u)}{\partial v^s}+}\nonumber\\
& & \mbox{\phantom{$\Phi^{ijk}$}} + \varphi^{is}(u,v)\frac{\partial \varphi^{jk}(v,w)}{\partial v^s}+\varphi^{si}(w,u)\frac{\partial \varphi^{kj}(w,v)}{\partial w^s}+\nonumber\\
& & \mbox{\phantom{$\Phi^{ijk}(u,v,w)=$}} + \varphi^{js}(v,w)\frac{\partial \varphi^{ki}(w,u)}{\partial w^s}+\varphi^{sj}(u,v)\frac{\partial \varphi^{ik}(u,w)}{\partial u^s}=0,\label{eqn19}
\eea
which supply triangular $r$-matrices for the Lie algebra $W_n$. We note that (\ref{eqn19}) is a system of $n^3$ equations for the $n^2$ unknowns formal series $\varphi^{ij}(u,v)\in{\Bbb R}[[u,v]]$, $i,j=1,\ldots,n$, with the property $\varphi^{ij}(u,v)=-\varphi^{ji}(v,u)$.
\begin{remark} We shall be interested first in the class of formal power series (polynomial) solutions $\varphi^{ij}(u,v)\in{\Bbb R}[[u,v]]$ of (\ref{eqn19}), for these will give rise to Poisson-Lie structures on $G_{n}$ and $G_{0n}$ (and Lie bialgebra structures on $W_n$ and its principal subalgebra). We shall also describe a more general encompassing class of formal Laurent series (rational) solutions $\varphi^{ij}(u,v)\in{\Bbb R}((u,v))$ of (\ref{eqn19}). 
However, the theorems formulated below hold vis-a-vis for solutions of  (\ref{eqn19}) sought in the class of smooth functions.  

Moreover, it follows from (\ref{eqn11}) that the group $G_{0n}$ acts on the space of solutions of the CYBE, $\Phi=0$. Naturally,  the problem of classification of solutions of (\ref{eqn19}) (within the set of  triangular $r$-matrices) is therefore equivalent to describing the set of $G_{0n}$-orbits (the moduli space of solutions of (\ref{eqn19})).
\end{remark}

The following theorem supplies a large class of solutions of (\ref{eqn19}). (Here, for clarity, we write explicitly the summation signs.)
\begin{theorem}
Let $F^i(u)\in{\Bbb R}((u^1,\ldots,u^n))$, $i=1,\ldots,n$, be $n$ arbitrary formal Laurent series (or more generally, components of a diffeomorphism $F:{{\Bbb R}^n}\to{\Bbb R}^n$).
Let $F_{*u}=\left(\frac{\partial F^i(u)}{\partial u^j}\right)$ be the
matrix of their partial derivatives (the derivative of the map F) and $F_{*u}^{-1}$ its inverse. Then the following $n^2$ formal (Laurent) series (smooth functions)
give a solution to the CYBE (\ref{eqn19}):
\bea
\varphi^{ij}(u,v)&=&\sum_{k=1}^{n}\sum_{l=1}^{n}\left(F_{*u}^{-1}\right)_{k}^{i}\left(F_{*v}^{-1}\right)_{l}^{j}\left[F^{k}\left(u\right)-F^{l}\left(v\right)\right]\label{eqn20}\\
&=&\sum_{k=1}^{n}\left(F_{*u}^{-1}\right)_{k}^{i}F^{k}\left(u\right)\sum_{l=1}^{n}\left(F_{*v}^{-1}\right)_{l}^{j}-\sum_{k=1}^{n}\left(F_{*u}^{-1}\right)_{k}^{i}\sum_{l=1}^{n}\left(F_{*v}^{-1}\right)_{l}^{j}F^{l}\left(v\right).\nonumber
\eea
\end{theorem}
\pf\ It is obvious that $\varphi^{ij}(u,v)=-\varphi^{ji}(v,u)$. The rest of the proof consists of a direct verification while using the following two facts. The first one is the identity
\be
\frac{\partial }{\partial u^s}\left(F_{*u}^{-1}\right)^{i}_{p}=-\left(F_{*u}^{-1}\right)^{i}_{q}\frac{\partial^{2} F^q(u)}{\partial u^s\partial u^j}\left(F_{*u}^{-1}\right)^{j}_{p}.\label{eqn20a}
\ee
The second obvious fact is that the quantities
\be
A_{pq}^i(u)=\frac{\partial^{2} F^i(u)}{\partial u^k\partial u^j}\left(F_{*u}^{-1}\right)^{k}_{p}\left(F_{*u}^{-1}\right)^{j}_{q}\label{eqn20b}
\ee
comprise the elements of a trivalent matrix which is symmetric with
respect to its lower  indices, i.e., $A_{pq}^i(u)=A_{qp}^i(u)$. 
Thus, after substitution of (\ref{eqn20}) into (\ref{eqn19}), use of
(\ref{eqn20a}), and regrouping of terms, we obtain for the left hand side  of (\ref{eqn19}):
\bea
%%\lefteqn{\left(F_{*u}^{-1}\right)_{k}^{s}\left(F_{*v}^{-1}\right)_{l}^{t}\left(F_{*w}^{-1}\right)_{n}^{p}[F^s(u)-F^p(w)+F^s(u)-F^p(w)+F^p(w)-F%%^t(v)+F^p(w)-}\nonumber\\
%%&&\mbox{\phantom{$(b-d)n_1+(c-a)n_2(b-d)n_1+$}} -F^t(v)+F^t(v)-F^s(u)+F^t(v)-F^s(u)]-\nonumber\\
\lefteqn{\left(F_{*u}^{-1}\right)^{k}_{r}\left(F_{*v}^{-1}\right)^{l}_{t}\left(F_{*w}^{-1}\right)^{n}_{s}\left\{\frac{\partial^{2} F^r(u)}{\partial u^m\partial u^j}\left(F_{*u}^{-1}\right)^{m}_{p}\left(F_{*u}^{-1}\right)^{j}_{q}\right\}\times}\nonumber\\
&&\times\{[F^p(u)-F^s(w)][F^q(u)-F^t(v)]+[F^t(v)-F^p(u)][F^q(u)-F^s(w)]\}+\nonumber\\
\lefteqn{+\left(F_{*u}^{-1}\right)^{k}_{s}\left(F_{*v}^{-1}\right)^{l}_{t}\left(F_{*w}^{-1}\right)^{n}_{r}\left\{\frac{\partial^{2} F^r(w)}{\partial w^m\partial w^j}\left(F_{*w}^{-1}\right)^{m}_{p}\left(F_{*w}^{-1}\right)^{j}_{q}\right\}\times}\nonumber\\
&& \times\{[F^s(u)-F^p(w)][F^q(w)-F^t(v)]+[F^p(w)-F^t(v)][F^q(w)-F^s(u)]\}+\nonumber\\
\lefteqn{+\left(F_{*u}^{-1}\right)^{k}_{t}\left(F_{*v}^{-1}\right)^{l}_{r}\left(F_{*w}^{-1}\right)^{n}_{s}\left\{\frac{\partial^{2} F^r(v)}{\partial v^m\partial v^j}\left(F_{*v}^{-1}\right)^{m}_{p}\left(F_{*v}^{-1}\right)^{j}_{q}\right\}\times}\nonumber\\
&& \times\{[F^s(w)-F^p(v)][F^q(v)-F^t(u)]+[F^p(v)-F^t(u)][F^q(v)-F^s(w)]\}\stackrel{\rm ?}{=}0.\nonumber
\eea
Using (\ref{eqn20b}) we rewrite the above as
\bea
\lefteqn{\left(F_{*u}^{-1}\right)^{k}_{r}\left(F_{*v}^{-1}\right)^{l}_{t}\left(F_{*w}^{-1}\right)^{n}_{p}A_{pq}^r(u){\mathcal F}^{pqst}_1(u,v,w)+}\nonumber\\
& & \mbox{\phantom{$\Phi^{ijk}$}} + \left(F_{*u}^{-1}\right)^{k}_{s}\left(F_{*v}^{-1}\right)^{l}_{t}\left(F_{*w}^{-1}\right)^{n}_{r}A_{pq}^r(w){\mathcal F}^{pqst}_2(u,v,w)+\nonumber\\
& & \mbox{\phantom{$\Phi^{ijk}(u,v,w)=$}} + \left(F_{*u}^{-1}\right)^{k}_{t}\left(F_{*v}^{-1}\right)^{l}_{r}\left(F_{*w}^{-1}\right)^{n}_{s}A_{pq}^r(v){\mathcal F}^{pqst}_3(u,v,w)\stackrel{\rm ?}{=}0,\label{eqn20ab}
\eea
where
\bea
{\mathcal F}^{pqst}_1(u,v,w)&=&[F^p(u)-F^s(w)][F^q(u)-F^t(v)]+[F^t(v)-F^p(u)][F^q(u)-F^s(w)],\nonumber\\
{\mathcal F}^{pqst}_2(u,v,w)&=&[F^s(u)-F^p(w)][F^q(w)-F^t(v)]+[F^p(w)-F^t(v)][F^q(w)-F^s(u)],\nonumber\\
{\mathcal F}^{pqst}_3(u,v,w)&=&[F^s(w)-F^p(v)][F^q(v)-F^t(u)]+[F^p(v)-F^t(u)][F^q(v)-F^s(w)].\nonumber
\eea
Now, it is easy to see that ${\mathcal F}^{pqst}_i(u,v,w)=-{\mathcal F}^{qpst}_i(u,v,w)$, $i=1,2,3$.
Using the fact that  $A_{pq}^r=A_{qp}^r$, we have 
\bea
A_{pq}^r(\cdot){\mathcal F}^{pqst}_i(u,v,w)&=&-A_{pq}^r(\cdot){\mathcal F}^{qpst}_i(u,v,w)\nonumber\\
&=&-A_{qp}^r(\cdot){\mathcal F}^{qpst}_i(u,v,w)\nonumber\\
&=&-A_{pq}^r(\cdot){\mathcal F}^{pqst}_i(u,v,w),\qquad i=1,2,3;\quad \cdot = u,w,v,\nonumber
\eea
from which it follows that $A_{pq}^r(\cdot){\mathcal F}^{pqst}_i(u,v,w)=0$, for $i=1,2,3$.
Thus each of the three
terms in (\ref{eqn20ab}) above is identically equal to zero. This establishes that (\ref{eqn20}) is a solution of the CYBE (\ref{eqn19}) and therefore completes the proof.\hfill \fbox{}

It turns out that the result of Theorem 3.1 is a particular case  of a more general statement.
\begin{theorem}
Let $\Psi^{i}(u)\in{\Bbb R}((u^1,\ldots,u^n))$, $i=1,\ldots,n$, be $n$ arbitrary given formal Laurent series (or a set of $n$ smooth functions), and let $\alpha,\beta\in{\Bbb R}$ be two arbitrary real numbers. Then a sufficient condition for
\be
\varphi^{ij}(u,v)=\Theta^i(u)\Psi^j(v)-\Theta^j(v)\Psi^i(u)\label{eqn21}
\ee
to be a solution of the CYBE (\ref{eqn19}), where  $\Theta^i(u)\in{\Bbb R}((u^1,\ldots,u^n))$ are $n$ Laurent series (smooth functions), is that $\Theta^i(u)$ satisfy the following system of $n$ linear partial differential equations
\be
\Psi^s(u)\frac{\partial \Theta^i(u)}{\partial u^s}-\Theta^s(u)\frac{\partial \Psi^i(u)}{\partial u^s}=\alpha\Theta^i(u)+\beta\Psi^i(u),\quad i=1,\ldots,n.\label{eqn22}
\ee
\end{theorem}
\pf\  After substituting (\ref{eqn21}) into (\ref{eqn19}) we obtain
\bea
\Theta^{i}(v)\Psi^{j}(w)\left[\Psi^{s}(u)\frac{\partial \Theta^{k}}{\partial u^s}-\Theta^{s}(u)\frac{\partial \Psi^{k}}{\partial u^s}\right]&+&\Theta^{j}(w)\Psi^{i}(v)\left[\Theta^{s}(u)\frac{\partial \Psi^{k}}{\partial u^s}-\Psi^{s}(u)\frac{\partial \Theta^{k}}{\partial u^s}\right]+\nonumber\\
&+&\mbox{cyclic}\left((i,u),(j,v),(k,w)\right)\stackrel{\rm ?}{=}0.\nonumber
\eea
From here, after the use of (\ref{eqn22}), one immediately concludes that the above expression  becomes an identity. \hfill \fbox{}
\begin{conjecture} (the necessity part)
All solutions of the CYBE are, in fact, described by the above theorem.
\end{conjecture}

For example, the solution (\ref{eqn20}) corresponds to the choice of $\alpha=0$, $\beta=1$, and $\Psi^{i}(u)=\sum_{s=1}^{n}\left(F_{*u}^{-1}\right)_{s}^{i}$, where $F_{*u}=\left(\frac{\partial F^i(u)}{\partial u^j}\right)$ for given  $F^i(u)$, $i=1,\ldots,n$.

\section{Poisson structures on ${\Bbb R}^n$}

With the construction of the Poisson-Lie groups in the previous two sections it is natural to ask whether there are Poisson structures on ${\Bbb R}^n$ such that
the action $FDiff({\Bbb R^n})\times {\Bbb R}^n\to {\Bbb R}^n$ of these Poisson-Lie groups of diffeomorphisms, defined by $u\mapsto X(u)$, is Poisson. We recall that an
action of a Poisson-Lie group $G$ on a Poisson manifold $M$ is said to be Poisson if the action map $G\times M\to M$ is Poisson, i.e., if  it respects the Poisson
structures on the source and image, where the source $G\times M$ is assumed to be equipped with the natural product Poisson structure. 

Let the Poisson-Lie group $FDiff({\Bbb R}^n)$ be equipped with a Poisson structure $\Omega$. Let $\alpha^{ij}$ be a Poisson bi-vector field on ${\Bbb R}^n$. It is by definition compatible with
the action of the Poisson-Lie group $FDiff({\Bbb R}^n)$  if and only if
\be
\alpha^{ij}(X(u))=\Omega^{ij}(j^{\infty}_{u}X,j^{\infty}_{u}X)+\left(X_{*u}\right)^i_k\left(X_{*u}\right)^j_l\alpha^{kl}(u),\label{eqn-1-newsec}
\ee
which is the condition that the map $u\mapsto X(u)$ is Poisson.
\begin{proposition} If the group  $FDiff({\Bbb R}^n)$ is equipped with a Poisson structure of type (\ref{eqn10}) then  the action of $FDiff({\Bbb R^n})$ on ${\Bbb R}^n$ by $u\mapsto X(u)$ induces a
Poisson structure $\alpha^{ij}$ on ${\Bbb R}^n$ given by $\alpha^{ij}(u)=-\varphi^{ij}(u,u)$, where $\varphi^{ij}(u,v)$ is any solution of the CYBE (\ref{eqn19}). In particular, if the $r$-matrix $\varphi^{ij}(u,v)$ is of type (\ref{eqn20}) then
\be
\alpha^{ij}(u)=-\varphi^{ij}(u,u)=-\sum_{k=1}^{n}\sum_{l=1}^{n}\left(F_{*u}^{-1}\right)_{k}^{i}\left(F_{*u}^{-1}\right)_{l}^{j}\left[F^{k}\left(u\right)-F^{l}\left(u\right)\right].\label{eqn-2-newsec}
\ee
\end{proposition}
\pf Let us rewrite equation (\ref{eqn-1-newsec}) with $\Omega^{ij}(j^{\infty}_{u}X,j^{\infty}_{u}X)$ given by formula (\ref{eqn10}):
$$
\alpha^{ij}(X(u))=\left(X_{*u}\right)^i_k\left(X_{*u}\right)^j_l\varphi^{kl}(u,u)-\varphi^{ij}(X(u),X(u))+\left(X_{*u}\right)^i_k\left(X_{*u}\right)^j_l\alpha^{kl}(u).
$$
Now, it is obvious that $\alpha^{ij}(u)=-\varphi^{ij}(u,u)$ is a solution of the above equation. To show that $\alpha^{ij}(u)$ is a Poisson tensor we evaluate the left hand side of the CYBE (\ref{eqn19}) on the diagonal $u=v=w$, which gives
\bea
\lefteqn{\varphi^{ks}(u,u)\left[\partial_{1s}\varphi^{ij}(u,u)-\partial_{1s}\varphi^{ji}(u,u)\right]+\varphi^{is}(u,u)\left[\partial_{1s}\varphi^{jk}(u,u)-\partial_{1s}\varphi^{kj}(u,u)\right]+}\nonumber\\
& &\mbox{\phantom{$\Phi^{ijk}(u,v,w)=$}} + \varphi^{js}(u,u)\left[\partial_{1s}\varphi^{ki}(u,u)-\partial_{1s}\varphi^{ik}(u,u)\right]=0.\label{eqn-3-newsec}
\eea
Here and below the notation $\partial_{1s}$ and $\partial_{2s}$ is used for the derivatives with respect to  the first and second (group of) arguments. We note that $\partial_{1s}\varphi^{ji}(u,u)=-\partial_{2s}\varphi^{ij}(u,u)$. Therefore for $\alpha^{ij}(u)=-\varphi^{ij}(u,u)$ we have $\frac{\partial \alpha^{ij}(u)}{\partial u^s}=-(\partial_{1s}\varphi^{ij}(u,u)+\partial_{2s}\varphi^{ij}(u,u))$. We conclude that the equation (\ref{eqn-3-newsec}) is equivalent to the Jacobi identities for $\alpha^{ij}$:
\be
\alpha^{ks}(u)\frac{\partial \alpha^{ij}(u)}{\partial u^s}+\alpha^{is}(u)\frac{\partial \alpha^{jk}(u)}{\partial u^s}+\alpha^{js}(u)\frac{\partial \alpha^{ki}(u)}{\partial u^s}=0.
\ee
Thus, $\alpha^{ij}$ is a Poisson bi-vector field on ${\Bbb R}^n$. Finally, from Theorem 3.1 follows that $\alpha^{ij}(u)$ given by formula (\ref{eqn-2-newsec}) satisfies the Jacobi identities. The proof is completed.\epf

We obtain the first non-trivial Poisson structure compatible with the action of the Poisson-Lie group of diffeomorphisms on ${\Bbb R}^2$ (on ${\Bbb R}^1$ we have only the trivial one). In the Appendix we give examples of such structures on ${\Bbb R}^2$ and ${\Bbb R}^3$.

\section{Poisson structures on the space of jets $J^{\infty}({\Bbb R^m},{\Bbb R^n})$}

In this section we give the answer to the following natural question. Let us consider the manifolds $R^m$ and $R^n$ and the infinite-dimensional manifold $J^{\infty}({\Bbb R^m},{\Bbb R^n})$  of $\infty$-jets of smooth maps from ${\Bbb R}^m$ to ${\Bbb R}^n$. Let $\Omega_m$ and $\Omega_n$ be given Poisson-Lie structures on the groups $FDiff({\Bbb R}^m)$ and $FDiff({\Bbb R}^n)$ respectively. The group $FDiff({\Bbb R}^n)\times FDiff({\Bbb R}^m)$ acts on the space $J^{\infty}({\Bbb R^m},{\Bbb R^n})$ in the following way. For every $X\in FDiff({\Bbb R}^m)$, $Y\in FDiff({\Bbb R}^n)$ and $F\in J^{\infty}({\Bbb R^m},{\Bbb R^n})$ the map $F\mapsto Y\circ F\circ X^{-1}$ defines an action of $FDiff({\Bbb R}^n)\times FDiff({\Bbb R}^m)$ on $J^{\infty}({\Bbb R^m},{\Bbb R^n})$ (left-right action). Here $X^{-1}$ is the inverse of $X$, i.e. $X^i(X^{-1}(u))=(X^{-1})^i(X(u))=u^i$, $i=1,\ldots,m$, for every point $u\in {\Bbb R^m}$ with coordinates $(u^1,\ldots,u^m)$. Let us take the space $FDiff({\Bbb R}^n)\times FDiff({\Bbb R}^m)$ with its natural product Poisson structure. We ask: are there Poisson structures $\Pi$ on the infinite-dimensional manifold $J^{\infty}({\Bbb R^m},{\Bbb R^n})$ such that the above action  is Poisson? The answer is found with the help of the following  three propositions.
\begin{proposition}
Let $\Omega_m$, $\Omega_n$ and $\Pi$ be Poisson structures on the groups $FDiff({\Bbb R}^m)$, $FDiff({\Bbb R}^n)$, and the space $J^{\infty}({\Bbb R^m},{\Bbb R^n})$ respectively. Then for every $X\in FDiff({\Bbb R}^m)$, $Y\in FDiff({\Bbb R}^n)$ and $F\in J^{\infty}({\Bbb R^m},{\Bbb R^n})$ the map $F\mapsto \overline{F}=Y\circ F\circ X^{-1}$ is Poisson if $\Pi$ satisfies
\bea
\lefteqn{\Pi^{ij}\left(j^{\infty}_{u}\overline{F},j^{\infty}_{v}\overline{F}\right)=\left(Y_{*F\circ X^{-1}(u)}\right)^i_k\left(Y_{*F\circ X^{-1}(v)}\right)^j_l\Pi^{kl}\left(j^{\infty}_{X^{-1}(u)}F,j^{\infty}_{X^{-1}(v)}F\right)+}\nonumber\\
& &\mbox{\phantom{$\Phi=$}} +\left(Y_{*F\circ X^{-1}(u)}\right)^i_k\left(F_{*X^{-1}(u)}\right)^k_{\sigma}\left(Y_{*F\circ X^{-1}(v)}\right)^j_l\left(F_{*X^{-1}(v)}\right)^l_{\rho}\Omega_m^{\sigma\rho}\left(j^{\infty}_{u}X^{-1},j^{\infty}_{v}X^{-1}\right)+\nonumber\\
& &\mbox{\phantom{$\Phi=$}} +\Omega_n^{ij}\left(j^{\infty}_{F\circ X^{-1}(u)}Y,j^{\infty}_{F\circ X^{-1}(v)}Y\right).\label{eqn-4-newsec}
\eea
\end{proposition}
We assume the convention that indices denoted with Greek letters run over the set $\{1,\ldots,m\}$ and those denoted with Latin letters run over the set $\{1,\ldots,n\}$.

The proof of Proposition 5.1 follows the same line as the proof of Proposition 3.1. Therefore we do not repeat the arguments. We proceed with the description of a class of solutions of the $n^2$ functional equations (\ref{eqn-4-newsec}) induced on the groups $FDiff({\Bbb R}^m)$ and $FDiff({\Bbb R}^n)$ by the Poisson structures $\Omega_m$ and $\Omega_n$ of type (\ref{eqn10}). 
\begin{proposition}
Let the Poisson-Lie structures on $FDiff({\Bbb R}^m)$ and $FDiff({\Bbb R}^n)$ be of type (\ref{eqn10}) and given by the formulae
\be
\Omega_m^{\sigma\rho}\left(j^{\infty}_{u}X,j^{\infty}_{v}X\right)=\left(X_{*u}\right)^{\sigma}_{\mu}\left(X_{*v}\right)^{\rho}_{\nu}\varphi_m^{\mu\nu}(u,v)-\varphi_m^{\sigma\rho}\left(X(u),X(v)\right),\label{eqn-5-newsec}
\ee
and
\be
\Omega_n^{ij}\left(j^{\infty}_{a}Y,j^{\infty}_{b}Y\right)=\left(Y_{*a}\right)^{i}_{k}\left(Y_{*b}\right)^{j}_{l}\varphi_n^{kl}(a,b)-\varphi_n^{ij}\left(Y(a),Y(b)\right),\label{eqn-6-newsec}
\ee
respectively, where $u,v\in {\Bbb R}^m$, $a,b\in {\Bbb R}^n$ and $\varphi_m^{\mu\nu}(u,v)$, $\varphi_n^{kl}(a,b)$ are the corresponding $r$-matrices. Then the equation (\ref{eqn-4-newsec}) has a solution of the form
\be 
\Pi^{ij}\left(j^{\infty}_{u}F,j^{\infty}_{v}F\right)=\left(F_{*u}\right)^{i}_{\mu}\left(F_{*v}\right)^{j}_{\nu}\varphi_m^{\mu\nu}(u,v)-\varphi_n^{ij}\left(F(u),F(v)\right).\label{eqn-7-newsec}
\ee
\end{proposition}
\pf
To show that (\ref{eqn-7-newsec}) is a solution of equation (\ref{eqn-4-newsec}) we substitute the expressions for $\Pi$, $\Omega_m$ and $\Omega_n$ in both of its left hand and right hand sides. After the substitution the left hand side reads
\be
\left(\overline{F}_{*u}\right)^{i}_{\mu}\left(\overline{F}_{*v}\right)^{j}_{\nu}\varphi_m^{\mu\nu}(u,v)-\varphi_n^{ij}\left(\overline{F}(u),\overline{F}(v)\right)\stackrel{\rm ?}{=}\label{eqn-8-newsec}
\ee
while the right hand side reads
\bea
\lefteqn{\stackrel{\rm ?}{=}\left(Y_{*F\circ X^{-1}(u)}\right)^i_k\left(Y_{*F\circ X^{-1}(v)}\right)^j_l\left(F_{*X^{-1}(u)}\right)^{k}_{\mu}\left(F_{*X^{-1}(v)}\right)^{l}_{\nu}\varphi_m^{\mu\nu}(u,v)-}\nonumber\\
& &\mbox{\phantom{$\Phi$}} -\left(Y_{*F\circ X^{-1}(u)}\right)^i_k\left(Y_{*F\circ X^{-1}(v)}\right)^j_l\varphi_n^{kl}\left(F\circ X^{-1}(u),F\circ X^{-1}(v)\right)+ \nonumber\\
& &\mbox{\phantom{$\Phi$}} + \left(Y_{*F\circ X^{-1}(u)}\right)^i_k\left(F_{*X^{-1}(u)})^k_{\sigma}(Y_{*F\circ X^{-1}(v)}\right)^j_l\left(F_{*X^{-1}(v)}\right)^l_{\rho}\left(X^{-1}_{*u}\right)^{\sigma}_{\mu}\left(X^{-1}_{*v}\right)^{\rho}_{\nu}\varphi_m^{\mu\nu}(u,v)-  \nonumber\\
& &\mbox{\phantom{$\Phi$}} -\left(Y_{*F\circ X^{-1}(u)}\right)^i_k\left(F_{*X^{-1}(u)}\right)^k_{\sigma}\left(Y_{*F\circ X^{-1}(v)}\right)^j_l\left(F_{*X^{-1}(v)}\right)^l_{\rho}\varphi_m^{\sigma\rho}\left(X^{-1}(u),X^{-1}(v)\right) +  \nonumber\\
& &\mbox{\phantom{$\Phi$}} + \left(Y_{*F\circ X^{-1}(u)}\right)^{i}_{k}\left(Y_{*F\circ X^{-1}(v)}\right)^{j}_{l}\varphi_n^{kl}\left(F\circ X^{-1}(u),F\circ X^{-1}(v)\right)-\nonumber\\
& &\mbox{\phantom{$\Phi$}} -\varphi_n^{ij}\left(Y\circ F\circ X^{-1}(u),Y\circ F\circ X^{-1}(v)\right)\nonumber \\
\lefteqn{=\left(Y_{*F\circ X^{-1}(u)})^i_k(F_{*X^{-1}(u)}\right)^k_{\sigma}\left(Y_{*F\circ X^{-1}(v)}\right)^j_l\left(F_{*X^{-1}(v)}\right)^l_{\rho}\left(X^{-1}_{*u}\right)^{\sigma}_{\mu}\left(X^{-1}_{*v}\right)^{\rho}_{\nu}\varphi_m^{\mu\nu}(u,v)-}\nonumber\\
& &\mbox{\phantom{$\Phi$}} -\varphi_n^{ij}\left(Y\circ F\circ X^{-1}(u),Y\circ F\circ X^{-1}(v)\right).\label{eqn-9-newsec}
\eea
Since we  have $\overline{F}=Y\circ F\circ X^{-1}$, the  chain rule implies $\overline{F}_{*u}=Y_{*F\circ X^{-1}(u)}F_{*X^{-1}(u)}X^{-1}_{*u}$. After inspecting the terms in (\ref{eqn-8-newsec}) and (\ref{eqn-9-newsec}) we see that we indeed obtain an identity. \epf

What remains is to find the sufficient condition for $\Pi$ to satisfy the Jacobi identity, i.e., to be Poisson.
We use again the representation $\left\{F^i(u),F^j(v)\right\}\equiv \Pi^{ij}\left(j^{\infty}_{u}F,j^{\infty}_{v}F\right)$. The Jacobi identity
then reads
\be
\{\{F^i(u),F^j(v)\},F^k(w)\}+\{\{F^k(w),F^i(u)\},F^j(v)\}+\{\{F^j(v),F^k(w)\},F^i(u)\}=0.\label{eqn-10-newsec}
\ee
This allows us to state and prove the following result.
\begin{proposition}
Let $\Pi$ be given by formula (\ref{eqn-7-newsec}). Then the Jacobi identity (\ref{eqn-10-newsec}) is equivalent to 
\be
\Phi_n^{ijk}(F(u),F(v),F(w))=(F_{*u})^{i}_{\mu}(F_{*v})^{j}_{\nu}(F_{*w})^{k}_{\sigma}\Phi_m^{\mu\nu\sigma}(u,v,w).\label{eqn-11-newsec}
\ee
That is, $\Pi$ given by (\ref{eqn-7-newsec}) is Poisson if and only if the invariant (\ref{eqn12}), $\Phi_n\in \mbox{Inv}_{G_{0n}}(W_n\widehat{\otimes}W_n\widehat{\otimes}W_n)$, is a  push-forward of the invariant $\Phi_m\in \mbox{Inv}_{G_{0m}}(W_m\widehat{\otimes}W_m\widehat{\otimes}W_m)$ along the map $F:{\Bbb R}^m\to{\Bbb R}^n$.
\end{proposition}
\pf
The proof consists of a calculation analogous to the calculation
involved in the proof of Proposition 3.3. We present only its essential part. Using formula (\ref{eqn-7-newsec}) we compute
$$
\{\{F^i(u),F^j(v)\},F^k(w)\}=\varphi_m^{\mu\nu}(u,v)\left\{\frac{\partial F^i(u)}{\partial u^{\mu}}\frac{\partial F^j(v)}{\partial v^{\nu}},F^k(w)\right\}-\left\{\varphi_n^{ij}(F(u),F(v)),F^k(w)\right\}
$$
$$
=\varphi_m^{\mu\nu}(u,v)\frac{\partial F^i(u)}{\partial u^{\mu}}\frac{\partial }{\partial v^{\nu}}\left\{F^j(v),F^k(w)\right\}+\varphi_m^{\mu\nu}(u,v)\frac{\partial F^j(v)}{\partial v^{\nu}}\frac{\partial }{\partial u^{\mu}}\left\{F^i(u),F^k(w)\right\}
$$
$$
-\partial_{1s}\varphi_n^{ij}(F(u),F(v))\left\{F^s(u),F^k(w)\right\}-\partial_{2s}\varphi_n^{ij}(F(u),F(v))\left\{F^s(v),F^k(w)\right\}
$$
$$
=\varphi_m^{\mu\nu}(u,v)\frac{\partial F^i(u)}{\partial u^{\mu}}\left[\frac{\partial^2F^j(v)}{\partial v^{\nu}\partial v^{\rho}}\frac{\partial F^k(w)}{\partial w^{\sigma}}\varphi_m^{\rho\sigma}(v,w)+\frac{\partial F^j(v)}{\partial v^{\rho}}\frac{\partial F^k(w)}{\partial w^{\sigma}}\partial_{1\nu}\varphi_m^{\rho\sigma}(v,w)\right]+
$$
$$
+\varphi_m^{\mu\nu}(u,v)\frac{\partial F^j(v)}{\partial v^{\nu}}\left[\frac{\partial^2F^i(u)}{\partial u^{\mu}\partial u^{\rho}}\frac{\partial F^k(w)}{\partial w^{\sigma}}\varphi_m^{\rho\sigma}(v,w)+\frac{\partial F^i(u)}{\partial u^{\rho}}\frac{\partial F^k(w)}{\partial w^{\sigma}}\partial_{1\mu}\varphi_m^{\rho\sigma}(u,w)\right]+
$$
$$
-\varphi_m^{\mu\nu}(u,v)\frac{\partial F^i(u)}{\partial u^{\mu}}\frac{\partial F^q(v)}{\partial v^{\nu}}\partial_{1q}\varphi_n^{jk}(F(v),F(w))-\varphi_m^{\mu\nu}(u,v)\frac{\partial F^j(v)}{\partial v^{\nu}}\frac{\partial F^q(u)}{\partial u^{\mu}}\partial_{1q}\varphi_n^{ik}(F(u),F(w))
$$
$$
-\partial_{1s}\varphi_n^{ij}(F(u),F(v))\frac{\partial F^s(u)}{\partial u^{\rho}}\frac{\partial F^k(w)}{\partial w^{\sigma}}\varphi_m^{\rho\sigma}(u,w)+\partial_{1s}\varphi_n^{ij}(F(u),F(v))\varphi_n^{sk}(F(u),F(w))
$$
$$
-\partial_{2s}\varphi_n^{ij}(F(u),F(v))\frac{\partial F^s(v)}{\partial v^{\rho}}\frac{\partial F^k(w)}{\partial w^{\sigma}}\varphi_m^{\rho\sigma}(u,w)+\partial_{2s}\varphi_n^{ij}(F(u),F(v))\varphi_n^{sk}(F(v),F(w)).
$$
After cyclic permutation of $(i,u)$, $(j,v)$ and $(k,w)$ we obtain analogous expressions for the remaining two terms in (\ref{eqn-10-newsec}). Adding the so obtained expressions for the three summands in (\ref{eqn-10-newsec}) we finally arrive at (\ref{eqn-11-newsec}) (after a short calculation consisting of matching terms which cancel each other out).\epf

We are now ready to make the following conclusion. Let $\varphi_m^{\mu\nu}(u,v)$ be a solution of the CYBE $\Phi_m^{\mu\nu\sigma}(u,v,w)=0$ and let $\varphi_n^{ij}(u,v)$ be a solution of the CYBE $\Phi_n^{ijk}(u,v,w)=0$. Then (\ref{eqn-11-newsec}) is identically satisfied for these two
$r$-matrices. Therefore each pair of such $r$-matrices determines a Poisson structure on the space $J^{\infty}({\Bbb R}^m,{\Bbb R}^n)$ of infinite jets, which is given by formula (\ref{eqn-7-newsec}).

\section{Classes of solutions of the CYBE}

We  now set out to describe (partially) the equivalence classes into which the solutions of the CYBE of type (\ref{eqn20}) split, modulo the action of $G_{0n}$. We start with the following theorem.
\begin{theorem}
Let ${\mathcal D}_{+}=\{D\in\mbox{Mat}(n,{\Bbb Z}_{+})\mid \mbox{det}D\neq 0\}$ be the set of non-singular $n\times n$ matrices with non-negative integer entries. Then to each $D\in{\mathcal D}_{+}$ corresponds a unique $G_{0n}$-orbit of solutions of the CYBE of type (\ref{eqn20}). Moreover, for each $D\in{\mathcal D}_{+}$ one can choose  a canonical triangular $r$-matrix $\varphi^{ij}_D(u,v)\in{\Bbb R}[u,v]$ representing the $G_{0n}$-orbit  labeled by $D$.
\end{theorem} 
\pf Let
\be
D=\pmatrix{d_{11}&d_{12}&\ldots&d_{1n}\cr d_{21}&d_{22}&\ldots&d_{2n}\cr \vdots&\vdots&\ddots&\vdots\cr d_{n1}&d_{n2}&\ldots&d_{nn}}\in \mbox{Mat}(n,{\Bbb Z}_{+}),\quad \mbox{det}D\neq 0,\label{eqn33}
\ee
be an non-singular $n\times n$ matrix with non-negative integer entries, $d_{ij}\in{\Bbb Z}_{+}$. To each such matrix $D$ there correspond the  $n$ formal Laurent series:
\bea
F^1(u)&=&(u^1)^{-d_{11}}(u^2)^{-d_{12}}\ldots(u^n)^{-d_{1n}}+\sum_{|{\bf k}|=-(d_{11}+\cdots +d_{1n})+1}^{\infty}f^1_{\bf k}{\bf u}^{\bf k},\nonumber\\
F^2(u)&=&(u^1)^{-d_{21}}(u^2)^{-d_{22}}\ldots(u^n)^{-d_{2n}}+\sum_{|{\bf k}|=-(d_{21}+\cdots +d_{2n})+1}^{\infty}f^2_{\bf k}{\bf u}^{\bf k},\nonumber\\
&\vdots&\label{eqn34}\\
F^n(u)&=&(u^1)^{-d_{n1}}(u^2)^{-d_{n2}}\ldots(u^n)^{-d_{nn}}+\sum_{|{\bf k}|=-(d_{n1}+\cdots +d_{nn})+1}^{\infty}f^n_{\bf k}{\bf u}^{\bf k},\nonumber
\eea
where ${\bf k}$ is a multi-index, ${\bf u}^{\bf k}=(u^1)^{k_{1}}(u^2)^{k_{2}}\ldots(u^n)^{k_{n}}$, $|{\bf k}|=k_1+\cdots +k_n$ and $f^i_{\bf k}\in{\Bbb R}$, $i=1,\ldots,n$, ${\bf k}\in{\Bbb Z}^n_{+}$, are arbitrary. The map so described gives a bijection between the set ${\mathcal D}_{+}$ and a set of
$G_{0n}$-orbits of solutions of the CYBE (\ref{eqn19}) of type (\ref{eqn20}). Indeed, let the $n$ formal Laurent series (\ref{eqn34}) be given and let  $X\in G_{0n}$ be a formal diffeomorphism such that
\bea
F^1(u)&=&(X^1(u))^{-d_{11}}(X^2(u))^{-d_{12}}\ldots(X^n(u))^{-d_{1n}},\nonumber\\
F^2(u)&=&(X^1(u))^{-d_{21}}(X^2(u))^{-d_{22}}\ldots(X^n(u))^{-d_{2n}},\nonumber\\
&\vdots&\label{eqn34a}\\
F^n(u)&=&(X^1(u))^{-d_{n1}}(X^2(u))^{-d_{n2}}\ldots(X^n(u))^{-d_{nn}}.\nonumber
\eea
Such a formal diffeomorphism exists and it is unique, granted  that $\mbox{det}D\neq 0$. If $\varphi^{ij}_D(u,v)$ is a solution of the CYBE (\ref{eqn19}) that is obtained from (\ref{eqn34a}) then it lies on the same $G_{0n}$-orbit of solutions as does the solution $\tilde{\varphi}^{ij}_D(u,v)={\varphi}^{ij}_D(X^{-1}(u),X^{-1}(v))$ which is obtained from 
\bea
\tilde{F}^1(u)&=&(u^1)^{-d_{11}}(u^2)^{-d_{12}}\ldots(u^n)^{-d_{1n}},\nonumber\\
\tilde{F}^2(u)&=&(u^1)^{-d_{21}}(u^2)^{-d_{22}}\ldots(u^n)^{-d_{2n}},\nonumber\\
&\vdots&\label{eqn34b}\\
\tilde{F}^n(u)&=&(u^1)^{-d_{n1}}(u^2)^{-d_{n2}}\ldots(u^n)^{-d_{nn}}.\nonumber
\eea
The solution $\tilde{\varphi}^{ij}_D(u,v)$ is the one that we choose to represent the orbit labeled by $D$.\epf

The general structure of the entries of the representative $r$-matrices described by Theorem 6.1 is given by the polynomials
\be
\varphi^{ij}(u,v)=u^iv^j\left[A^{ij}_1{\bf u}^{{\bf d}_1}+\cdots+A^{ij}_n{\bf u}^{{\bf d}_n}-A^{ji}_1{\bf v}^{{\bf d}_1}-\cdots-A^{ji}_n{\bf v}^{{\bf d}_n}\right],\quad i,j=1,\ldots,n,\label{eqn34bc}
\ee
where ${\bf d}_i=(d_{i1},d_{i2},\ldots,d_{in})$ are the rows of the matrix $D$. The coefficients $A^{ij}_n$ are expressed  in terms of the $(n-1)\times(n-1)$ minors of the matrix $D$. In particular, the
Poisson tensor $\alpha^{ij}$ on ${\Bbb R}^n$ induced by these $r$-matrices has polynomial  components and is given by
\be
\alpha^{ij}(u)=u^iu^j\left[\left(A^{ij}_1-A^{ji}_1\right){\bf u}^{{\bf d}_1}+\cdots+\left(A^{ij}_n-A^{ji}_n\right){\bf u}^{{\bf d}_n}\right],\quad i,j=1,\ldots,n.\label{eqn34bcc}
\ee
In the Appendix we give explicit formulae for $r$-matrices of type (\ref{eqn34bc}) and the corresponding Poisson-Lie brackets for the groups $FDiff({\Bbb R}^n)$, $n=1,2,3$. 

We term {\em non-singular} the isomorphism classes of solutions corresponding to $F^i(u)\in{\Bbb R}((u^1,\ldots,u^n))$, $i=1,\ldots,n$ as in (\ref{eqn34}). Let ${\mathcal D}_{+}^{sing}=\{D\in\mbox{Mat}(n,{\Bbb Z}_{+})\mid \mbox{det}D=0\}$ be the set of singular $n\times n$ matrices with non-negative integer entries. We term {\em degenerate} the set  of $r$-matrices obtained from (\ref{eqn34bc}) by a choice of $D\in{\mathcal D}_{+}^{sing}$. In this case there
exist integers $k_1,k_2,\ldots,k_n\in{\Bbb Z}$ such that
\be
k_1{\bf d}_1+k_2{\bf d}_2+\cdots+k_n{\bf d}_n=0.
\ee
Analyzing the structure of the set of degenerate solutions should be interesting.

There is one more $G_{0n}$-orbit of solutions of the CYBE of type (\ref{eqn20}) that can be canonically represented by polynomial $r$-matrices. It is described as follows. 
\begin{theorem}
To the set $\left\{A\in\mbox{Mat}(n,{\Bbb R})\mid\mbox{det}A\neq 0\right\}$ of non-singular $n\times n$ real matrices, there corresponds a unique $G_{0n}$-orbit of solutions of the CYBE of type (\ref{eqn20}). Moreover, there is a canonical triangular $r$-matrix $\tilde{\varphi}^{ij}(u,v)\in{\Bbb R}[u,v]$ representing this $G_{0n}$-orbit.
\end{theorem}
\pf
For  any non-singular matrix
\be
A=\pmatrix{a_{11}&a_{12}&\ldots&a_{1n}\cr a_{21}&a_{22}&\ldots&a_{2n}\cr \vdots&\vdots&\ddots&\vdots\cr a_{n1}&a_{n2}&\ldots&a_{nn}}\in \mbox{Mat}(n,{\Bbb R}),\quad \mbox{det}A\neq 0,\label{eqn33ac}
\ee
let the $n$ formal power series 
\bea
F^1(u)&=&a_{11}u^1+a_{12}u^2+\cdots +a_{1n}u^n+\sum_{|{\bf k}|>1}^{\infty}f^1_{\bf k}{\bf u}^{\bf k},\nonumber\\
F^2(u)&=&a_{21}u^1+a_{22}u^2+\cdots +a_{2n}u^n+\sum_{|{\bf k}|>1}^{\infty}f^2_{\bf k}{\bf u}^{\bf k},\nonumber\\
&\vdots&\label{eqn34aa}\\
F^n(u)&=&a_{n1}u^1+a_{n2}u^2+\cdots +a_{nn}u^n+\sum_{|{\bf k}|>1}^{\infty}f^n_{\bf k}{\bf u}^{\bf k},\nonumber
\eea
be given. Then the non-singularity of $A$ implies that there is a unique formal diffeomorphism $X\in G_{0n}$ such that
\bea
F^1(u)&=&X^{1}(u),\nonumber\\
F^2(u)&=&X^{2}(u),\nonumber\\
&\vdots&\label{eqn34ab}\\
F^n(u)&=&X^{n}(u).\nonumber
\eea
Thus a solution $\varphi^{ij}(u,v)$ of the CYBE (\ref{eqn19}) that is obtained from (\ref{eqn34aa}) can be transformed to the solution $\tilde{\varphi}^{ij}(u,v)={\varphi}^{ij}(X^{-1}(u),X^{-1}(v))$ on the same $G_{0n}$-orbit that  is obtained from
\bea
\tilde{F}^1(u)&=&u^{1},\nonumber\\
\tilde{F}^2(u)&=&u^{2},\nonumber\\
&\vdots&\label{eqn34ba}\\
\tilde{F}^n(u)&=&u^{n}.\nonumber
\eea
Moreover, the representative solution has a particularly simple form given by
\be
\tilde{\varphi}^{ij}(u,v)=u^i-v^j,\quad i,j=1,2,\ldots,n.\label{eqnsigma}
\ee\epf

Clearly the $G_{0n}$-orbit represented by (\ref{eqnsigma}) can be adjoint to the set of orbits described by Theorem 6.1 by adjoining the matrix $-I$ to the set ${\mathcal D}_{+}$, where $I$ is the  $n\times n$ identity matrix.

So far we considered $G_{0n}$-orbits of solutions $\varphi^{ij}(u,v)$ of the CYBE which are formal power series (polynomials). These power series (polynomial) solutions give rise to Poisson-Lie structures on the group of diffeomorphisms $FDiff({\Bbb R}^n)$. The set of their $G_{0n}$-orbits is in one-to-one correspondence with the set ${\mathcal D}_{+}\cup\{-I\}$. If we remove the restriction that the $\varphi^{ij}(u,v)$ should be formal power series then we enter a wider class of solutions given by formal Laurent series. In general, the latter can be obtained as follows.

Let ${\mathcal D}=\{D\in\mbox{Mat}(n,{\Bbb Z})\mid \mbox{det}D\neq 0\}$ be the set of non-singular $n\times n$ matrices with integer entries. Then to each matrix 
\be
D=\pmatrix{d_{11}&d_{12}&\ldots&d_{1n}\cr d_{21}&d_{22}&\ldots&d_{2n}\cr \vdots&\vdots&\ddots&\vdots\cr d_{n1}&d_{n2}&\ldots&d_{nn}}\in \mbox{Mat}(n,{\Bbb Z}),\quad \mbox{det}D\neq 0,\label{eqn33}
\ee
we can associate the totality ${\mathcal F}_{D}(n)$ of $n$-tuples of formal (Laurent) series 
\bea
F^1(u)&=&\sum_{k_1=-d_{11}}^{\infty}\sum_{k_2=-d_{12}}^{\infty}\cdots\sum_{k_n=-d_{1n}}^{\infty}f^1_{k_1k_2\ldots k_n}(u^1)^{k_{1}}(u^2)^{k_{2}}\ldots(u^n)^{k_{n}},\nonumber\\
F^2(u)&=&\sum_{k_1=-d_{21}}^{\infty}\sum_{k_2=-d_{22}}^{\infty}\cdots\sum_{k_n=-d_{2n}}^{\infty}f^2_{k_1k_2\ldots k_n}(u^1)^{k_{1}}(u^2)^{k_{2}}\ldots(u^n)^{k_{n}},\nonumber\\
&\vdots&\label{eqn-totality}\\
F^n(u)&=&\sum_{k_1=-d_{n1}}^{\infty}\sum_{k_2=-d_{n2}}^{\infty}\cdots\sum_{k_n=-d_{nn}}^{\infty}f^n_{k_1k_2\ldots k_n}(u^1)^{k_{1}}(u^2)^{k_{2}}\ldots(u^n)^{k_{n}}.\nonumber
\eea
Thus, if $D\not\in{\mathcal D}_{+}\cup\{-I\}$, any element of ${\mathcal F}_{D}(n)$ gives rise to a Laurent series solution of (\ref{eqn19}) obtained from formula (\ref{eqn20}). For instance,
they give rise to classes of Lie bialgebra structures  on the Lie algebra Vect$({\Bbb T}^n)$ of vector fields on the torus ${\Bbb T}^n$. We shall address the classification of these in a separate paper.

\section{Classification of all Poisson-Lie structures on $FDiff({\Bbb R}^1)$}

In this section we consider in more detail the case of $G_{1}=FDiff({\Bbb R}^1)$ and its Lie algebra $W_1$. This is the simplest case. Here one is able to advance  on the solution of the classification problem further than in the general case of $FDiff({\Bbb R}^n)$, and $W_n$, for $n>1$. In particular, there exists an effective
method of constructing all solutions of the CYBE and describing their moduli space. The system of equations (\ref{eqn11}) reduces to the single equation 
\be
\Phi(X(u),X(v),X(w))=X'(u)X'(v)X'(w)\Phi(u,v,w)\label{eqn23}
\ee
in this case. Here $X(u)\in{\Bbb R}[[u]]$ will be assumed within this section to be a formal diffeomorphism $X\in G_{01}=FDiff_0({\Bbb R}^1)$,
\bea
\lefteqn{\Phi(w,u,v)=\varphi(u,v)\bigl[ {\partial}_u\varphi(w,u)+{\partial}_v\varphi(w,v)\bigr]+}\nonumber\\
& & \mbox{\phantom{$\Phi^{ijk}\Phi^{ijk}\Phi^{ijk}\Phi^{ijk}$}} +\varphi(v,w)\bigl[{\partial}_v\varphi(u,v)+{\partial}_w\varphi(u,w)\bigr]+\nonumber\\
& & \mbox{\phantom{$\Phi^{ijk}(u,v,w)=\Phi^{ijk}(u,v,w)=$}} +\varphi(w,u)\bigl[{\partial}_w\varphi(v,w)+{\partial}_u\varphi(v,u)\bigr],\label{eqn24}
\eea
and we assume that $\varphi(u,v)\in{\Bbb R}[[u,v]]$. 
It is not difficult to show \cite{Stoyanov:paper_JMP} that $\mbox{Inv}_{G_{01}}(W_1\widehat{\otimes}W_1\widehat{\otimes}W_1)=0$.  Therefore all Poisson-Lie structures on $G_{1}$ (and all Lie bialgebra structures on $W_1$) are obtained from  solutions of the CYBE, $\Phi=0$. Moreover, the group $G_{01}$ acts on the space of solutions, as follows from (\ref{eqn23}). To find all solutions of the CYBE and the space of moduli we proceed as follows. For $w=v+\epsilon$ let us
expand $\Phi(u,v,v+\epsilon)$ near the diagonal\footnote{The author was convinced of the effectiveness of this technique (in many situations similar to the one at hand) by P.Etingof, in a private communication.}
$$
\Phi(u,v,v+\epsilon)=\frac {\partial\Phi}{\partial w}{\left|\right.}_{w=v}\epsilon+\frac {1}{2}\frac {\partial^2\Phi}{\partial w^2}{\mid}_{w=v}\epsilon^2+\cdots.
$$
We therefore obtain that $\Phi(u,v,w)=0$ is equivalent to the infinite system of differential equations $\frac {\partial\Phi}{\partial w}{\mid}_{w=v}=0,\frac {\partial^2\Phi}{
\partial w^2}{\mid}_{w=v}=0,\ldots$. However, to find all solutions of $\Phi(u,v,w)=0$ it is enough to only find all solutions
of the (weaker) equation $\frac {\partial\Phi}{\partial w}{\mid}_{w=v}=0$ and then show that they all are also solutions of $\Phi(u,v,w)=0$. We have thus,
\be
\frac {\partial\Phi}{\partial w}{\mid}_{w=v}=\partial_v\varphi(u,v)\partial_u\varphi(u,v)-\partial_u\partial_v\varphi(u,v)+f'(v)\varphi(u,v)-f(v)\partial_v\varphi(u,v)=0,\label{eqn25}
\ee
where $f(v):=\partial_1\varphi(v,v)$, $f'(v)=\partial_1^2\varphi(v,v)$ and we have used the fact that $\partial_1\partial_2\varphi(v,v)=0$. Here $\partial_1$ and $\partial_2$ denote the derivatives with respect to the first and second arguments. Away from the diagonal $u=v$ the last equation is then equivalent to
\be
\partial_v\left[-\frac{\partial_u\varphi(u,v)}{\varphi(u,v)}+\frac{f(v)}{\varphi(u,v)}\right]=0,\label{eqn26}
\ee
from where it follows that $\varphi(u,v)$ satisfies the linear partial differential equation
\be
\partial_u\varphi(u,v)+G(u)\varphi(u,v)=f(v),\label{eqn27}
\ee
where $G(u)\in{\Bbb R}[[u,u^{-1}]]$ is an arbitrary formal Laurent series. The general solution of the last equation is
\be
\varphi(u,v)=\left[f(v)\int_{0}^{u}\exp({{\int_{0}^{t}}G(s)ds})dt+H(v)\right]\exp({-\int_{0}^{u}G(s)ds}),\label{eqn28}
\ee
where $H(v)\in{\Bbb R}[[u]]$ is an arbitrary power (Taylor) series. Next, the condition $\varphi(u,u)=0$ implies that $H(u)=-f(u)\int_{0}^{u}\exp({{\int_{0}^{t}G(s)ds})dt}$, and further, the skew-symmetry  $\varphi(u,v)=-\varphi(v,u)$ implies that $G(u)=-\frac{f'(u)}{f(u)}$. Therefore we obtain the following formula for the general solution of (\ref{eqn26}):
\be
\varphi(u,v)=f(u)f(v)\left[\int_{0}^{u}\frac{dt}{f(t)}-\int_{0}^{v}\frac{dt}{f(t)}\right].\label{eqn29}
\ee
Now, let $F(u)\in{\Bbb R}((u))$ be a formal Laurent series such that $f(u)=1/F'(u)$. Then
\be
\varphi(u,v)=\frac{1}{F'(u)}\frac{1}{F'(v)}\left[F(u)-F(v)\right].\label{eqn30}
\ee
It is now straightforward to verify that the so obtained general solution of (\ref{eqn26}) satisfies the CYBE (\ref{eqn24}), and thus gives {\it its} general solution. Write
\be
F(u)=\sum_{n=-d}^{\infty}C_{n}u^n=C_{-d}\frac {1}{u^{d}}+C_{-d+1}\frac {1}{u^{d-1}}+\cdots,\label{eqn31}
\ee
for some $d\in{\Bbb Z}$. However, from the choice of $F(u)$ follows that $d\ge -1$, since $f(u)\in{\Bbb R}[[u]]$ is a formal series. Let us choose  a formal
diffeomorphism $X(u)$ such that $F(u)=-(X(u))^{-d}$. Under the inverse of this diffeomorphism  the solution (\ref{eqn30}) gets transformed to 
\be
\tilde{\varphi}(u,v)=\frac {1}{d^2}u^{d+1}v^{d+1}(v^{-d}-u^{-d}).\label{eqn32}
\ee
Note that $d=0$ corresponds to the trivial solution. We therefore conclude that elements of the moduli space of all solutions of the CYBE for $W_1$ are in one to one correspondence with the set ${\Bbb Z}_{+}\cup\{-1\}$. Finally, it was shown in \cite{Stoyanov:paper_JMP} that $H^1(W_1,W_1\widehat{\otimes}W_1)=0$, which implies that all Lie bialgebra structures on $W_1$ are coboundary and come from an $r$-matrix. As a result we obtain a classification of all Lie bialgebra structures on $W_1$ and equivalently all Poisson-Lie structures on $FDiff({\Bbb R^1})$.
\begin{remark}
Removing the restriction that $f(u)$ above is a formal power series
allows us to consider arbitrary formal Laurent series $F(u)\in{\Bbb R}((u))$ or formal power series $F(u)\in{\Bbb R}[[u]]$ and thus obtain
$r$-matrices $\varphi(u,v)\in{\Bbb R}((u,v))$ that are Laurent series. In particular, we obtain countably many isomorphism classes of Lie bialgebra
structures on the algebra Vect$({\Bbb S^1})$ of vector fields on
${\Bbb S^1}$, each labeled by an element  $d\in{\Bbb Z}$ of the set of integers, and represented by (\ref{eqn32}).
\end{remark}

\appendix
\section{Appendix}

In this appendix we describe several examples. We give the explicit formulae for the representative
classes of Poisson brackets on the
space of $\infty$-jets of formal diffeomorphisms of ${\Bbb R}^1$ as
well as ${\Bbb R}^2$, which are obtained from the  class (\ref{eqn20}) of 
solutions   of the CYBE (\ref{eqn19}). We also present the representative $r$-matrices for the  Lie algebra $W_3$ of formal vector fields on ${\Bbb R}^3$ and the class of Poisson brackets they induce on ${\Bbb R}^3$. The general formulae for the Poisson
bracket on the space of $\infty$-jets of formal diffeomorphisms of ${\Bbb R}^n$ are rather cumbersome and we refrain from writing them here.

We start with $FDiff({\Bbb R}^1)$. According to the results obtained in the previous section to each $d\in{\Bbb Z}_{+}\cup\{-1\}$ corresponds an isomorphism
class of solutions of the CYBE $\Phi=0$ from (\ref{eqn24}), and represented by (\ref{eqn32}). Each of these classes gives rise to a class of Poisson-Lie structures on $FDiff({\Bbb R}^1)$ given by
\bea
\lefteqn{\left\{x_i,x_j\right\}_d=(i-d)jx_jx_{i-d}-i(j-d)x_ix_{j-d}+}\nonumber\\
&&\mbox{\phantom{$\{x_i,x_j\}$}} +x_i\sum_{(\stackrel{\mbox{\tiny all partitions}}{s_1+\cdots+s_{d+1}=j})}x_{s_1}\ldots x_{s_{d+1}}-x_j\sum_{(\stackrel{\mbox{\tiny all partitions}}{s_1+\cdots+s_{d+1}=i})}x_{s_1}\ldots x_{s_{d+1}},\nonumber
\eea
where $i,j$ roam over the set ${\Bbb Z}_{+}$. They have been described in earlier publications \cite{KS1:paper,Stoyanov:paper_JMP}. In the formula above we implicitly
assume that $x_i=0$, whenever $i<0$.

Consider now the case of $FDiff({\Bbb R}^2)$.
Let
$$
D=\pmatrix{a&b\cr c&d}\in \mbox{Mat}(2,{\Bbb Z}_{+}),\quad \mbox{det}D\neq 0.
$$
Let us associate to this  matrix the Laurent monomials
\bea
F^1(u)&=&(u^1)^{-a}(u^2)^{-b}\nonumber\\
F^2(u)&=&(u^1)^{-c}(u^2)^{-d}.\nonumber
\eea
Then, corresponding to these monomials, we obtain  the following $r$-matrix of type (\ref{eqn20}):
\bea
\varphi^{11}(u,v)&=&u^1v^1\left\{(b-d)d\left[(u^1)^a(u^2)^b-(v^1)^a(v^2)^b\right]-(b-d)b\left[(u^1)^c(u^2)^d-(v^1)^c(v^2)^d\right]\right\}\nonumber\\
\varphi^{12}(u,v)&=&u^1v^2\left\{(b-d)\left[c(v^1)^a(v^2)^b-a(v^1)^c(v^2)^d\right]+(c-a)\left[d(u^1)^a(u^2)^b-b(u^1)^c(u^2)^d\right]\right\}\nonumber\\
\varphi^{21}(u,v)&=&u^2v^1\left\{(c-a)\left[-d(v^1)^a(v^2)^b+b(v^1)^c(v^2)^d\right]+(b-d)\left[-c(u^1)^a(u^2)^b+a(u^1)^c(u^2)^d\right]\right\}\nonumber\\
\varphi^{22}(u,v)&=&u^2v^2\left\{(a-c)c\left[(u^1)^a(u^2)^b-(v^1)^a(v^2)^b\right]-(a-c)a\left[(u^1)^c(u^2)^d-(v^1)^c(v^2)^d\right]\right\}.\nonumber
\eea

For a germ of a diffeomorphism $X:{\Bbb R}^2\to {\Bbb R}^2$,  given by
$X^i(u)=\sum^{\infty}_{m_1=0}\sum^{\infty}_{m_2=0}x^{i}_{m_1m_2}(u^1)^{m_1}(u^2)^{m_2}$, 
$i=1,2$ the formula for the Poisson structures of type (\ref{eqn10}) and corresponding
to the above $r$-matrix imply the following formulae for the Poisson
brackets of the coordinates of $X$ as an element of the space $J^{\infty}({\Bbb R}^2,{\Bbb R}^2)$:
\bea
\lefteqn{\left\{x^1_{m_1m_2},x^1_{n_1n_2}\right\}=\left[(b-d)n_1+(c-a)n_2\right][\left(dm_1-cm_2-\mbox{det}D\right)x^{1}_{m_1-am_2-b}x^{1}_{n_1n_2}+}\nonumber\\\nonumber\\
&&\mbox{\phantom{$(b-d)n_1+(c-a)n_2(b-d)n_1+$}} +\left(-bm_1+am_2-\mbox{det}D\right)x^{1}_{m_1-cm_2-d}x^{1}_{n_1n_2}]+\nonumber\\\nonumber\\
&&\mbox{\phantom{$-a)n_2(b-d)$}} +\left[(b-d)m_1+(c-a)m_2\right][\left(-dn_1+cn_2+\mbox{det}D\right)x^{1}_{m_1m_2}x^{1}_{n_1-an_2-b}+\nonumber\\\nonumber\\
&&\mbox{\phantom{$-a)n(b-d)n_1+(c-a)n_2(b-d)$}} +\left(bn_1-an_2+\mbox{det}D\right)x^{1}_{m_1m_2}x^{1}_{n_1-cn_2-d}]+\nonumber\\\nonumber\\
&&+d(b-d)x^1_{n_1n_2}\sum_{(\stackrel{\mbox{\tiny all partitions}}{s_1+\cdots+q_b=m_1})}\sum_{(\stackrel{\mbox{\tiny all partitions}}{p_1+\cdots+r_b=m_2})}x^1_{s_1p_1}\ldots x^1_{s_{a+1}p_{a+1}}x^2_{q_1r_1}\ldots x^2_{q_br_b}-\nonumber\\
&&-d(b-d)x^1_{m_1m_2}\sum_{(\stackrel{\mbox{\tiny all partitions}}{s_1+\cdots+q_b=n_1})}\sum_{(\stackrel{\mbox{\tiny all partitions}}{p_1+\cdots+r_b=n_2})}x^1_{s_1p_1}\ldots x^1_{s_{a+1}p_{a+1}}x^2_{q_1r_1}\ldots x^2_{q_br_b}-\nonumber\\
&&-b(b-d)x^1_{n_1n_2}\sum_{(\stackrel{\mbox{\tiny all partitions}}{s_1+\cdots+q_d=m_1})}\sum_{(\stackrel{\mbox{\tiny all partitions}}{p_1+\cdots+r_d=m_2})}x^1_{s_1p_1}\ldots x^1_{s_{c+1}p_{c+1}}x^2_{q_1r_1}\ldots x^2_{q_dr_d}+\nonumber\\
&&+b(b-d)x^1_{m_1m_2}\sum_{(\stackrel{\mbox{\tiny all partitions}}{s_1+\cdots+q_d=n_1})}\sum_{(\stackrel{\mbox{\tiny all partitions}}{p_1+\cdots+r_d=n_2})}x^1_{s_1p_1}\ldots x^1_{s_{c+1}p_{c+1}}x^2_{q_1r_1}\ldots x^2_{q_dr_d},\nonumber
\eea

\bea
\lefteqn{\left\{x^1_{m_1m_2},x^2_{n_1n_2}\right\}=\left[(b-d)n_1+(c-a)n_2\right][\left(dm_1-cm_2-\mbox{det}D\right)x^{1}_{m_1-am_2-b}x^{2}_{n_1n_2}+}\nonumber\\\nonumber\\
&&\mbox{\phantom{$(b-d)n_1+(c-a)n_2(b-d)$}} +\left(-bm_1+am_2-\mbox{det}D\right)x^{1}_{m_1-cm_2-d}x^{2}_{n_1n_2}]+\nonumber\\\nonumber\\
&&\mbox{\phantom{$-a)n_2(b-d)$}} +\left[(b-d)m_1+(c-a)m_2\right][\left(-dn_1+cn_2+\mbox{det}D\right)x^{1}_{m_1m_2}x^{2}_{n_1-an_2-b}+\nonumber\\\nonumber\\
&&\mbox{\phantom{$-a)n(b-d)n_1+(c-a)n_2(b-d)$}} +\left(bn_1-an_2+\mbox{det}D\right)x^{1}_{m_1m_2}x^{2}_{n_1-cn_2-d}]+\nonumber\\\nonumber\\
&&+c(b-d)x^1_{m_1m_2}\sum_{(\stackrel{\mbox{\tiny all partitions}}{s_1+\cdots+q_{b+1}=n_1})}\sum_{(\stackrel{\mbox{\tiny all partitions}}{p_1+\cdots+r_{b+1}=n_2})}x^1_{s_1p_1}\ldots x^1_{s_ap_a}x^2_{q_1r_1}\ldots x^2_{q_{b+1}r_{b+1}}-\nonumber\\
&&-a(b-d)x^1_{m_1m_2}\sum_{(\stackrel{\mbox{\tiny all partitions}}{s_1+\cdots+q_{d+1}=n_1})}\sum_{(\stackrel{\mbox{\tiny all partitions}}{p_1+\cdots+r_{d+1}=n_2})}x^1_{s_1p_1}\ldots x^1_{s_cp_c}x^2_{q_1r_1}\ldots x^2_{q_{d+1}r_{d+1}}+\nonumber\\
&&+d(c-a)x^2_{n_1n_2}\sum_{(\stackrel{\mbox{\tiny all partitions}}{s_1+\cdots+q_b=m_1})}\sum_{(\stackrel{\mbox{\tiny all partitions}}{p_1+\cdots+r_b=m_2})}x^1_{s_1p_1}\ldots x^1_{s_{a+1}p_{a+1}}x^2_{q_1r_1}\ldots x^2_{q_br_b}-\nonumber\\
&&-b(c-a)x^2_{n_1n_2}\sum_{(\stackrel{\mbox{\tiny all partitions}}{s_1+\cdots+q_d=m_1})}\sum_{(\stackrel{\mbox{\tiny all partitions}}{p_1+\cdots+r_d=m_2})}x^1_{s_1p_1}\ldots x^1_{s_{c+1}p_{c+1}}x^2_{q_1r_1}\ldots x^2_{q_dr_d},\nonumber
\eea

\bea
\lefteqn{\left\{x^2_{m_1m_2},x^2_{n_1n_2}\right\}=\left[(b-d)n_1+(c-a)n_2\right][\left(dm_1-cm_2-\mbox{det}D\right)x^{2}_{m_1-am_2-b}x^{2}_{n_1n_2}+}\nonumber\\\nonumber\\
&&\mbox{\phantom{$(b-d)n_1+(c-a)n_2(b-d)$}} +\left(-bm_1+am_2-\mbox{det}D\right)x^{2}_{m_1-cm_2-d}x^{2}_{n_1n_2}]+\nonumber\\\nonumber\\
&&\mbox{\phantom{$-a)n_2(b-d)$}} +\left[(b-d)m_1+(c-a)m_2\right][\left(-dn_1+cn_2+\mbox{det}D\right)x^{2}_{m_1m_2}x^{2}_{n_1-an_2-b}+\nonumber\\\nonumber\\
&&\mbox{\phantom{$-a)n(b-d)n_1+(c-a)n_2(b-d)$}} +\left(bn_1-an_2+\mbox{det}D\right)x^{2}_{m_1m_2}x^{2}_{n_1-cn_2-d}]+\nonumber\\\nonumber\\
&&+c(a-c)x^2_{n_1n_2}\sum_{(\stackrel{\mbox{\tiny all partitions}}{s_1+\cdots+q_{b+1}=m_1})}\sum_{(\stackrel{\mbox{\tiny all partitions}}{p_1+\cdots+r_{b+1}=m_2})}x^1_{s_1p_1}\ldots x^1_{s_ap_a}x^2_{q_1r_1}\ldots x^2_{q_{b+1}r_{b+1}}-\nonumber\\
&&-c(a-c)x^2_{m_1m_2}\sum_{(\stackrel{\mbox{\tiny all partitions}}{s_1+\cdots+q_{b+1}=n_1})}\sum_{(\stackrel{\mbox{\tiny all partitions}}{p_1+\cdots+r_{b+1}=n_2})}x^1_{s_1p_1}\ldots x^1_{s_ap_a}x^2_{q_1r_1}\ldots x^2_{q_{b+1}r_{b+1}}-\nonumber\\
&&-a(c-a)x^2_{n_1n_2}\sum_{(\stackrel{\mbox{\tiny all partitions}}{s_1+\cdots+q_{d+1}=m_1})}\sum_{(\stackrel{\mbox{\tiny all partitions}}{p_1+\cdots+r_{d+1}=m_2})}x^1_{s_1p_1}\ldots x^1_{s_{c}p_{c}}x^2_{q_1r_1}\ldots x^2_{q_{d+1}r_{d+1}}+\nonumber\\
&&+a(c-a)x^2_{m_1m_2}\sum_{(\stackrel{\mbox{\tiny all partitions}}{s_1+\cdots+q_{d+1}=n_1})}\sum_{(\stackrel{\mbox{\tiny all partitions}}{p_1+\cdots+r_{d+1}=n_2})}x^1_{s_1p_1}\ldots x^1_{s_{c}p_{c}}x^2_{q_1r_1}\ldots x^2_{q_{d+1}r_{d+1}}.\nonumber
\eea
In the formulae above we implicitly assume that $x^i_{i_1i_2}=0$ if either $i_1<0$ or $i_2<0$. We also note that these formulae introduce a structure of a Poisson algebra on the infinitely generated commutative associative algebra ${\mathcal A}={\Bbb R}[[{X^{coor}}]]$, where ${X^{coor}}=\left\{x^i_{i_1i_2}\mid (i_1,i_2)\in{\Bbb Z}^2_{+},i=1,2\right\}$.

The Poisson bracket (compatible with the action of $X$) induced on ${\Bbb R}^2$ from the $r$-matrix 
corresponding to the integer matrix $D$ is 
\be
\left\{u^1,u^2\right\}=-\varphi^{12}(u,u)=u^1u^2\left[\left(u^1\right)^{d}\left(u^2\right)^{c}-\left(u^1\right)^{b}\left(u^2\right)^{a}\right].\label{appx-3}
\ee

Our next example is a polynomial $r$-matrix of type (\ref{eqn20}) for the Lie algebra $W_3$ of formal vector fields on ${\Bbb R}^3$. Chose a $3\times 3$
matrix 
$$
D=\pmatrix{d_{11}&d_{12}&d_{13}\cr d_{21}&d_{22}&d_{23}\cr d_{31}&d_{32}&d_{33}}\in \mbox{Mat}(3,{\Bbb Z}_{+}),\quad \mbox{det}D\neq 0,
$$
with non-negative integer entries. Let us associate to this matrix
the Laurent monomials:
\bea
F^1(u)&=&\left(u^1\right)^{-d_{11}}\left(u^2\right)^{-d_{12}}\left(u^3\right)^{-d_{13}},\nonumber\\
F^2(u)&=&\left(u^1\right)^{-d_{21}}\left(u^2\right)^{-d_{22}}\left(u^3\right)^{-d_{23}},\nonumber\\
F^3(u)&=&\left(u^1\right)^{-d_{31}}\left(u^2\right)^{-d_{32}}\left(u^3\right)^{-d_{33}}.\nonumber
\eea
To the above choices corresponds the following $r$-matrix of type
(\ref{eqn20}):
\bea
\varphi^{11}(u,v)&=&D_{11}(D_{11}-D_{21}+D_{31})u^1v^1\left[\frac{1}{F^1(v)}-\frac{1}{F^1(u)}\right]-\nonumber\\
&&\mbox{\phantom{$(b-d)$}} -D_{21}(D_{11}-D_{21}+D_{31})u^1v^1\left[\frac{1}{F^2(v)}-\frac{1}{F^2(u)}\right]+\nonumber\\
&&\mbox{\phantom{$(b-d)$}} +D_{31}(D_{11}-D_{21}+D_{31})u^1v^1\left[\frac{1}{F^3(v)}-\frac{1}{F^3(u)}\right],\nonumber\\
\varphi^{22}(u,v)&=&-D_{12}(-D_{12}+D_{22}-D_{32})u^2v^2\left[\frac{1}{F^1(v)}-\frac{1}{F^1(u)}\right]+\nonumber\\
&&\mbox{\phantom{$(b-d)$}} +D_{22}(-D_{12}+D_{22}-D_{32})u^2v^2\left[\frac{1}{F^2(v)}-\frac{1}{F^2(u)}\right]-\nonumber\\
&&\mbox{\phantom{$(b-d)$}} -D_{32}(-D_{12}+D_{22}-D_{32})u^2v^2\left[\frac{1}{F^3(v)}-\frac{1}{F^3(u)}\right],\nonumber\\
\varphi^{33}(u,v)&=&D_{13}(D_{13}-D_{23}+D_{33})u^3v^3\left[\frac{1}{F^1(v)}-\frac{1}{F^1(u)}\right]-\nonumber\\
&&\mbox{\phantom{$(b-d)$}} -D_{23}(D_{13}-D_{23}+D_{33})u^3v^3\left[\frac{1}{F^2(v)}-\frac{1}{F^2(u)}\right]+\nonumber\\
&&\mbox{\phantom{$(b-d)$}} +D_{33}(D_{13}-D_{23}+D_{33})u^3v^3\left[\frac{1}{F^3(v)}-\frac{1}{F^3(u)}\right],\nonumber\\
\varphi^{12}(u,v)&=&(D_{11}-D_{21}+D_{31})u^1v^2\left[-\frac{D_{12}}{F^{1}(v)}+\frac{D_{22}}{F^{2}(v)}-\frac{D_{32}}{F^{3}(v)}\right]+\nonumber\\
&&\mbox{\phantom{$(b-d)$}} +(D_{12}-D_{22}+D_{32})u^1v^2\left[\frac{D_{11}}{F^{1}(u)}-\frac{D_{21}}{F^{2}(u)}+\frac{D_{31}}{F^{3}(u)}\right],\nonumber\\
\varphi^{13}(u,v)&=&(D_{11}-D_{21}+D_{31})u^1v^3\left[\frac{D_{13}}{F^{1}(v)}-\frac{D_{23}}{F^{2}(v)}+\frac{D_{33}}{F^{3}(v)}\right]-\nonumber\\
&&\mbox{\phantom{$(b-d)$}} -(D_{13}-D_{23}+D_{33})u^1v^3\left[\frac{D_{11}}{F^{1}(u)}-\frac{D_{21}}{F^{2}(u)}+\frac{D_{31}}{F^{3}(u)}\right],\nonumber\\
\varphi^{23}(u,v)&=&(-D_{12}+D_{22}-D_{32})u^2v^3\left[\frac{D_{13}}{F^{1}(v)}-\frac{D_{23}}{F^{2}(v)}+\frac{D_{33}}{F^{3}(v)}\right]-\nonumber\\
&&\mbox{\phantom{$(b-d)$}} -(D_{13}-D_{23}+D_{33})u^2v^3\left[-\frac{D_{21}}{F^{1}(u)}+\frac{D_{22}}{F^{2}(u)}-\frac{D_{32}}{F^{3}(u)}\right].\nonumber
\eea
In the above formulae $D_{11}=d_{22}d_{33}-d_{23}d_{32}$,
$D_{12}=d_{21}d_{33}-d_{23}d_{31}$, etc.,  denote the $2\times 2$ minors of the matrix 
$D$. The induced Poisson brackets
$\left\{u^i,u^j\right\}=-\varphi^{ij}(u,u)$, $i,j=1,2,3$, on ${\Bbb R}^3$
are then described by the formulae:
\bea
\left\{u^{1},u^{2}\right\}&=&\left[D_{12}(D_{21}-D_{31})-D_{11}(D_{22}-D_{32})\right]\left(u^1\right)^{d_{11}+1}\left(u^2\right)^{d_{12}+1}\left(u^3\right)^{d_{13}}+\nonumber\\\nonumber\\
&&\mbox{\phantom{$()$}} +\left[D_{22}(D_{11}+D_{31})-D_{21}(D_{12}+D_{32})\right]\left(u^1\right)^{d_{21}+1}\left(u^2\right)^{d_{22}+1}\left(u^3\right)^{d_{23}}+\nonumber\\\nonumber\\
&&\mbox{\phantom{$()$}} +\left[D_{32}(-D_{11}+D_{21})+D_{31}(D_{12}-D_{22})\right]\left(u^1\right)^{d_{31}+1}\left(u^2\right)^{d_{32}+1}\left(u^3\right)^{d_{33}},\nonumber\\\nonumber\\
\left\{u^{1},u^{3}\right\}&=&\left[D_{13}(-D_{21}+D_{31})+D_{11}(D_{23}-D_{33})\right]\left(u^1\right)^{d_{11}+1}\left(u^2\right)^{d_{12}}\left(u^3\right)^{d_{13}+1}+\nonumber\\\nonumber\\
&&\mbox{\phantom{$()$}} +\left[-D_{23}(D_{11}+D_{31})+D_{21}(D_{13}+D_{33})\right]\left(u^1\right)^{d_{21}+1}\left(u^2\right)^{d_{22}}\left(u^3\right)^{d_{23}+1}+\nonumber\\\nonumber\\
&&\mbox{\phantom{$()$}} +\left[D_{33}(D_{11}-D_{21})+D_{31}(D_{23}-D_{13})\right]\left(u^1\right)^{d_{31}+1}\left(u^2\right)^{d_{32}}\left(u^3\right)^{d_{33}+1},\nonumber\\\nonumber\\
\left\{u^{2},u^{3}\right\}&=&\left[D_{13}(D_{22}-D_{32})+D_{12}(D_{33}-D_{23})\right]\left(u^1\right)^{d_{11}}\left(u^2\right)^{d_{12}+1}\left(u^3\right)^{d_{13}+1}+\nonumber\\\nonumber\\
&&\mbox{\phantom{$()$}} +\left[D_{23}(D_{12}+D_{32})-D_{22}(D_{13}+D_{33})\right]\left(u^1\right)^{d_{21}}\left(u^2\right)^{d_{22}+1}\left(u^3\right)^{d_{23}+1}+\nonumber\\\nonumber\\
&&\mbox{\phantom{$()$}} +\left[D_{33}(D_{22}-D_{12})+D_{32}(D_{13}-D_{23})\right]\left(u^1\right)^{d_{31}}\left(u^2\right)^{d_{32}+1}\left(u^3\right)^{d_{33}+1}.\nonumber
\eea
In particular, if $D$ is the diagonal matrix
$$
D=\pmatrix{a&0&0\cr 0&b&0\cr 0&0&c}\in \mbox{Mat}(3,{\Bbb Z}_{+}),\quad \mbox{det}D\neq 0,
$$
then we have
\bea
\left\{u^{1},u^{2}\right\}&=&cu^{1}u^{2}\left[\left(u^{2}\right)^{b}-\left(u^{1}\right)^{a}\right],\nonumber\\\nonumber\\
\left\{u^1{},u^{3}\right\}&=&bu^{1}u^{3}\left[\left(u^{3}\right)^{c}-\left(u^{1}\right)^{a}\right],\nonumber\\\nonumber\\
\left\{u^{2},u^{3}\right\}&=&au^{2}u^{3}\left[\left(u^{3}\right)^{c}-\left(u^{2}\right)^{b}\right].\nonumber
\eea

% (\stackrel{\mbox{\tiny all partitions}}{})

\label{kupershmidt_5-lp}

\end{document}